\documentclass[11pt]{article}
\usepackage{amsmath,amssymb,amsthm}
\usepackage{graphicx}
\usepackage{color}
\theoremstyle{plain}
\newtheorem{theo}{Theorem}[section]
\newtheorem*{theo*}{Theorem}
\newtheorem{lem}[theo]{Lemma}
\newtheorem*{lem*}{Lemma}

\newtheorem{cor}[theo]{Corollary}
\theoremstyle{definition}
\newtheorem{remark}[theo]{Remark}

\newtheorem*{remark*}{Remark}

\newtheorem*{ex*}{Example}
\theoremstyle{remark}
\theoremstyle{definition}

\newtheorem*{defi*}{Definition}
\newtheorem*{nota*}{Notation}
\numberwithin{equation}{section}
\newcommand{\mm}{\hspace{-0.5mm} }
\newcommand{\Linn}{\\ \hline
  &  &  &  &  & &  &   & & & &    \\    \vspace{-5mm}    \\
 \hline }
 1

 2

\usepackage{amsmath}
\usepackage{amssymb}

\newtheorem{theorem}{Theorem}[section]
\newtheorem{lemma}{Lemma}[section]

\usepackage{color}
\usepackage[normalem]{ulem}

\begin{document}

\title{Leading terms of velocity and its gradient of the stationary rotational viscous incompressible flows with
nonzero velocity at infinity.
}

\author {
Paul Deuring\footnote{Univ Lille
Nord de France, 59000 Lille, France; ULCO, LMPA, 62228 Calais c\'edex, France.}
,\
Stanislav Kra\v cmar\footnote{
Department of Technical Mathematics, Czech Technical University,
Karlovo n\'{a}m. 13, 121 35 Prague 2, Czech Republic}
,\  \v S\'arka Ne\v casov\' a\footnote{Mathematical Institute
of the Academy of Sciences of the Czech Republic, \v Zitn\' a~25,
115 67 Prague 1, Czech Republic}}
\date{ }
\maketitle

\vspace{2ex}
\begin{abstract}
We consider the Navier-Stokes system with Oseen and rotational terms describing the stationary flow of a viscous incompressible fluid around a rigid body moving at a constant velocity and rotating at a constant angular velocity. In a previous paper, we prove a representation formula for weak solutions of the system. Here the representation formula is used to get an asymptotic expansion of respectively velocity and its gradient, and  to establish pointwise decay estimates of remainder terms. Our results are based on a fundamental solution proposed by Guenther and Thomann \cite{GT}.
We thus present a different  approach to this result, besides the one, given by Kyed \cite{K}.



\end{abstract}



\section{Introduction}
The aim of this paper is to find the asymptotic structure, particularly the leading terms,
of the velocity part of the solution to the system
 \begin{equation}
\left.\begin{array}{rcl}
- \mu\Delta u (z) - (U+\omega \times z) \cdot  \nabla u (z) + \omega \times u (z) + u  \cdot \nabla u(z)
+\nabla \pi(z) &= &f(z), \\ \hbox{\rm div } u (z)& =& 0,
\label{nonlin} \end{array}\right\}
\end{equation}
\begin{eqnarray} \label{1.2}
u(x)\to 0\;\; \mbox{for}\;\; |x|\to \infty .
\end{eqnarray}
This system
describes the stationary flow of a viscous incompressible fluid around a rigid
body moving at a constant velocity and rotating at a constant angular velocity.
We refer to \cite{G2} for more details on the physical background of (\ref{nonlin}).
Here we only indicate that $ \mbox{$\mathfrak D$} \subset \mathbb{R}^3 $ is an
open bounded set describing the rigid body, the vector $ U \in \mathbb{R}^3 \backslash\{0\} $
represents the constant translational velocity of this body, the vector
$ \omega \in \mathbb{R}^3 \backslash\{0\} $ stands for its constant angular velocity,
and $\mu $ denotes the constant kinematic viscosity of the fluid. The
given function $f: \mathbb{R}^3 \backslash \overline{ \mbox{$\mathfrak D$} } \mapsto
\mathbb{R}^3 $ describes a body force, and the unknowns
$u: \mathbb{R}^3 \backslash \overline{ \mbox{$\mathfrak D$} } \mapsto
\mathbb{R}^3 $ and $\pi : \mathbb{R}^3 \backslash \overline{ \mbox{$\mathfrak D$} } \mapsto
\mathbb{R}$
correspond respectively
to the velocity and pressure field of the fluid.
We assume that $U \cdot \omega  \ne 0$. Then, according to \cite{GK1}, without loss
of generality we may replace (\ref{nonlin}) by the normalized system
\begin{eqnarray} \label{1.10}
L(u) + \tau (u \cdot \nabla )u + \nabla \pi = f,
\;\;
\mbox{div}\, u = 0
\quad \mbox{in}\;\; \mathbb{R}^3 \backslash \overline{ \mbox{$\mathfrak D$} },
\end{eqnarray}
where the differential operator $L$ is defined by
\begin{equation*} \label{1.21}
\begin{split} &
L(u)(z):= - \Delta u(z) + \tau \hspace{1pt}  \partial _1u(z) - ( \omega \times z) \cdot \nabla u(z)
+ \omega \times u(z)
\\ &
\mbox{for}\;\; u \in W^{2,1}_{loc}( U)^3,
\; z \in U,
\;
U \subset \mathbb{R}^3 \; \mbox{open},
\end{split}
\end{equation*}
with $\tau \in (0, \infty ) $ (Reynolds number) and $\omega = \varrho (1,0,0)$ for some
$\varrho \in \mathbb{R}\setminus\{0\} $  (Taylor number).

Suppose that $f \in L^{p_0}( \mathbb{R}^3 )^3$ for some $p_0 \in (1, \infty )$ and $f$ has
compact support. Further suppose there is a pair of functions $(u,\,\pi )$ with
$u \in L^6( \overline{ \mbox{$\mathfrak D$}} ^c)^3,\; \nabla u\in L^2( \overline{ \mbox{$\mathfrak D$}} ^c)^9$
and $\pi \in L^2_{loc}( \overline{ \mbox{$\mathfrak D$} }^c)$ satisfying (\ref{nonlin}) in the distributional
sense (''Leray solution''). Such a solution exists under suitable assumptions on
$\partial \mbox{$\mathfrak D$} ,\, u| \partial \mbox{$\mathfrak D$} $ and $p_0$
(\cite[Theorem XI.3.1]{Galdineu}).
Note that the condition
$u \in L^6( \overline{ \mbox{$\mathfrak D$} }^c)^3,\; \nabla u\in L^2( \overline{ \mbox{$\mathfrak D$}} ^c)^9$
means in particular that (\ref{1.2}) holds in a weak sense; compare \cite[Theorem II.5.1]{Ga1}.
In this situation, it was shown by Galdi and Kyed \cite{GK1} that
\begin{eqnarray} \label{1.4}
| \partial ^{\alpha }u(x)|=O \bigl[\, \bigl(\, |x| \, s_{\tau }(x) \,\bigr) ^{-1-| \alpha |/2}  \,\bigr]
\quad (|x|\to \infty ),
\end{eqnarray}
where $ \alpha \in \mathbb{N} ^3_0 $ with $ |\alpha |:=\alpha _1+\alpha _2+\alpha _3 \le 1$
(decay of $u$ and $\nabla u$).
The term $s_{\tau }(x) $ in (\ref{1.4}) is defined by
\begin{eqnarray} \label{1.4a}
s_{\tau }(x):= 1+\tau\,(|x|-x_1) \quad (x \in \mathbb{R}^3 ).
\end{eqnarray}
Its presence in (\ref{1.4}) may be considered as a mathematical manifestation of the wake extending
downstream behind the rigid body. Even in the linear nonrotational case, that is, in the case
of solutions to the Oseen system
\begin{eqnarray} \label{1.5}
-\Delta u + \tau \, \partial _1 u + \nabla \pi = f,
\;\;
\mbox{div}\, u =0,
\end{eqnarray}
the velocity cannot be expected to decay more rapidly than
$\bigl(\, |x| \, s_{\tau }(x)  \,\bigr) ^{-1}$ for $|x|\to \infty $, nor its gradient more rapidly
than $\bigl(\, |x| \, s_{\tau }(x) \,\bigr) ^{-3/2} $ (\cite{KNP}). Therefore the decay rate in (\ref{1.4})
should be best possible in the present case, too.
 By Kyed  \cite{K} it was shown that
\begin{eqnarray} \label{k}
u(x)= \mbox{$\mathfrak O$}  (x) \cdot \alpha +R(x), \quad
\nabla u(x)= \nabla \mbox{$\mathfrak O $}  (x) \cdot \alpha +S(x),
\end{eqnarray}
where $\mbox{$\mathfrak O$} $ is the fundamental solution of the stationary Oseen system,
$\alpha $ represents the force
\mbox{$\mathfrak F$} exerted by the liquid on the body, and $R$ and $S$ are some remainder terms decaying
faster than $ \mbox{$\mathfrak O\,\mbox{and} \,\nabla \mathfrak O,$}$ respectively, as $|x|\to \infty $.

In the work at hand, we also derive an asymptotic expansion of respectively $u$ and $\nabla u$. These
expansions -- stated in Theorem 3.1 below -- differ in two respects from those presented in \cite{K} and indicated
in (\ref{k}). Firstly, our leading term is less explicit than the term $\mbox{$\mathfrak O$}  \cdot \alpha $ in (\ref{k}).
Instead of the fundamental solution $\mbox{$\mathfrak O$}  $ of the \underline{stationary} Oseen system, we use
the time integral of the fundamental solution of the \underline{evolutionary} Oseen system multiplied by a rotation depending on time.
Secondly, and this is an aspect which goes beyond the theory in \cite{K}, we establish pointwise decay estimates
of our remainder terms (see ({\ref{T4.2.20aNL}})), whereas in \cite{K}, it is only shown that the function $R$ in (\ref{k}) belongs
to $L^p( \mathbb{R}^3 \backslash B_S)^3$ for $p \in (4/3,\, \infty ),$ and $S$ to $L^p( \mathbb{R}^3 \backslash
B_S)^9$ for $p \in (1, \infty )$, where $B_S$ is an open ball with sufficiently large radius $S>0$.
Interestingly, by integrating the decay rates in ({\ref{T4.2.20aNL}}) and using Lemma \ref{lemma1}, we find that our remainder terms belong to the same $L^p$-spaces.

We further indicate that our results are derived by an approach different from the one in \cite{K}: whereas
the theory in \cite{K} reduces (\ref{k}) to estimates of solutions to the time-periodic Oseen system in the
whole space $\mathbb{R}^3 $, our results are based on a representation formula of solutions to
({\ref{1.10}}) (see Theorem {\ref{theorem5.20}}).
As a consequence of our approach, our remainder terms are expressed explicitly in terms of $u,\,
\pi$ and $f$. In particular, sharpening (\ref{k}), we obtain that
$S= \nabla R$.

Our access is made difficult by the structure
of the Guenther-Thomann fundamental solution. In fact, as was already pointed out in
\cite{FHM} for the case $\tau =0$,
a fundamental solution $ \mbox{$\mathfrak Z$}(x,y)$ to (\ref{1.10}) cannot be bounded by $c \, |x-y| ^{-1} $
uniformly in $x,y \in \mathbb{R}^3 $ with $|x|$ and $|y|$ large, contrary to what may be expected
in view of the situation in the Stokes and Oseen case. Actually it seems that no uniform bound
$c \, |x-y|^{- \epsilon } $ exists, for whatever $ \epsilon \in (0, \infty ) $.

In  \cite{DKN1} -- \cite{DKN4}, we proved a representation formula, a decay estimate as in  (\ref{1.4}), and asymptotic expansions for
weak solutions of the linearized problem
\begin{eqnarray*}
L(u) +\nabla \pi =f, \quad \mbox{div}\, u =0, \quad u(x)\to \infty \;\;(|x|\to \infty ),
\end{eqnarray*}
as well as a
representation formula for weak solutions of (\ref{1.10}). In the context of these papers, a weak solution
$(u,\pi )$ of (\ref{1.10})  is
characterized by the assumptions that $u$ is $L^6$-integrable outside
a ball containing $ \overline{ \mbox{$\mathfrak D$} }$, and $ \nabla u$ and $\pi$ are $L^2$-integrable
outside such a ball. In \cite{DKN6}, we extended the
results from \cite{DKN1} -- \cite{DKN4}  from weak solutions to Leray solutions.
%

In \cite {DKN5}, we considered the nonlinear problem (\ref{1.10}), deriving optimal rates of decay as in (\ref{1.4}) for the velocity
and its gradient, on the basis of the representation formula proved in \cite{DKN2} and \cite{DKN6} and restated below as
Theorem \ref{theorem5.20}.


The asymptotic behavior of purely rotating case was studied by Farwig, Hishida; see  \cite{FH091}, \cite{FH092} for the linear case and  \cite{FH2, FGK} in the nonlinear one.

Concerning further articles related to the work at hand, we mention
\cite{AC},
\cite{Fa3}, \cite{FGNT}, \cite{FH},
\cite{FKN1} -- \cite{FN},
\cite{G3}, \cite{GK2}, \cite{GK3}, \cite{GHH},
\cite{H1} -- \cite{KNP-JAP},
\cite{KrPe1} -- \cite{K3},
\cite{Ne2}, \cite{NS}.

Let us briefly indicate how we will proceed in the following. In Section 2 we will present various auxiliary results. Section 3 deals with the main theorem - leading term for the velocity field and its  gradient.

\section{Notation and preliminaries}
The open bounded set  $ \mbox{$\mathfrak D$} \subset \mathbb{R}^3 $
introduced in Section 1 will be kept fixed throughout. We assume its
boundary
$ \partial \mbox{$\mathfrak D$} $ to be of class $C^2$,
and we denote its outward unit normal by
$n^{ ( \mbox{$\mathfrak D$}  )} $.
The numbers $ \tau $ and $\varrho $ and the vector $ \omega $ also
introduced in Section 1
will be kept fixed, too.
Define the matrix $ \Omega \in \mathbb{R}^{3 \times 3}  $ by
\begin{eqnarray*}&&
\Omega :=
\left( \begin{array}{rrr}
0 & - \omega_3  & \omega_2  \\
\omega _3 & 0 & - \omega _1\\
- \omega _2 & \omega _1 & 0
\end{array}
\right)
=
\varrho \hspace{1pt}  \left( \begin{array}{rrr}
0 & 0  & 0  \\
0& 0 & - 1\\
0 & 1 & 0
\end{array}
\right),
\end{eqnarray*}
so that $ \omega  \times x = \Omega \cdot x$ for $x \in \mathbb{R}^3
$.

Let us denote $s(x):=s_1(x)=1+(|x|-x_1).$ We recall that the function $s_{\tau }$
was defined in Section 1, as was the notation $| \alpha |$ for the
length of a multi-index $\alpha \in \mathbb{N} _0^3$.
If $ A \subset \mathbb{R}^3 $, we write $A^c$ for the complement $
\mathbb{R}^3 \backslash
A$ of $A$.
The open ball centered at $x \in \mathbb{R}^3 $ and with radius
$r>0$ is denoted by $B_r(x)$. If $x=0$, we will write $B_r$ instead of
$B_r(0)$. Put $ e_1:=(1,0,0)$.
Let $x \times y$ denote the usual vector product
of $x,y \in \mathbb{R}^3 $.
For $ T \in (0, \infty ) $, set $ \mbox{$\mathfrak D$} _T:=
B_T \backslash \overline{ \mbox{$\mathfrak D$} }$ (''truncated exterior
domain'').
By the symbol \mbox{$\mathfrak C$}, we denote constants only depending
on $ \mbox{$\mathfrak D$} ,\; \tau $ or
$ \omega $. We write $ \mbox{$\mathfrak C$} ( \beta _1,..., \beta _n)$
for positive constants
that additionally depend on parameters $ \beta _1,\, ...,\, \beta _n
\in \mathbb{R} $,
for some $ n \in \mathbb{N} $.
As usual, $C(\gamma_1,\dots,\gamma_n) $ means a positive constant only depending on $\gamma_1,\dots,\gamma_n.$



We will further use the
ensuing estimate, which was proved {in} \cite{Farwig2}.

\begin{lem}\label{lemma1}
Let $ \beta \in (1, \infty )$. Then
$
\int_{ \partial B_{r} } s_{ \tau }(x) ^{ - \beta }\; do_x \le
\mbox{$\mathfrak C$}(\beta )
\,  r
$ \
for
$ r \in (0, \infty ).$
\end{lem}

We begin by introducing the fundamental solutions used in what follows.
We set
\begin{equation*}
\begin{split}
&
K(x,t)= (4\pi t)^{-3/2} e^{-\frac{|x|^2}{4t}}, \ x\in \mathbb R^3, \ t\in (0,\infty),\label{eq5}
\\ &
N_{jk}(x) = x_jx_k|x|^{-2},\  x\in\mathbb R^3 \setminus
\{0\},
\\ &
\Lambda_{jk} (x,t) = K(x,t) \bigg (\delta_{jk}-N_{jk}(x) -
{}_1F_1\bigg(1,5/2,\frac{|x|^2}{4t}\bigg)
\,
\left(\delta_{jk}/3 - N_{jk}(x)\right)\bigg),
\\ &
\ x \in \mathbb R^3\setminus \{0\},\
t \in (0,\infty), \ j,k \in \{1,2,3\},
\\ &
{}_1F_1(1,c,u):= \sum_{n=0} ^{ \infty } \bigl(\, \Gamma (c)/ \Gamma (n+c) \,\bigr)  \cdot
u^n \quad \mbox{for}\;\; u \in \mathbb{R} ,\; c \in (0, \infty ),
\\ &
 \mbox{where }\Gamma \mbox { denotes the usual Gamma function. In the following, the letter
$\Gamma $ will stand for }\\ & \mbox{the matrix-valued function defined by}
\\ & (\Gamma_{jk} (y,z,t))_{1\le j,k\le 3} :=(\Lambda_{rs}(y-\tau \, t\,e_1-e^{-t\Omega}\cdot z,t))_{1\le r, s\le 3}
\cdot e^{-t\Omega}, 
\\ &
y,z \in \mathbb{R}^3 ,\ t \in (0, \infty ) \; \mbox{with}\;
y-\tau \, t\,e_1-e^{-t\Omega}\cdot z\ne 0.
\\ &
E_{4j}(x):= (4 \hspace{1pt}  \pi ) ^{-1} \hspace{1pt}  x_j
\hspace{1pt}  |x| ^{-3}, \quad 1\le j\le 3,\; x \in
\mathbb{R}^3 \backslash\{0\} .
\end{split}
\end{equation*}


Our following lemma restates \cite[Corollary 3.1]{DKN1}:

\begin{lem}
 The function $\Gamma $ may be continuously extended to a function from
$C ^{ \infty } \bigl(\,  \mathbb{R}^3 \times \mathbb{R}^3 \times  (0, \infty ) \,\bigr) .$
\end{lem}

According to \cite[Theorem 3.1]{DKN1}, we have
\begin{lem}\label{lemma2}
$
\int^{\infty}_0 |\Gamma_{jk} (y,z,t)| dt < \infty
$
for $y,z \in\mathbb R^3$ with $y\not=z$, $1\le j$, $k\le 3$.
\end{lem}

\medskip

Thus we may define
$$
\mbox{$\mathfrak Z$} _{jk} (y,z): = \int^{\infty}_0 \Gamma_{jk} (y,z,t) dt
$$
for $y,z \in\mathbb R^3$ with $y\not=z$, $1\le j$, $k \le 3$.
\medskip

The matrix-valued function \mbox{$\mathfrak Z$} constitutes the velocity part of the fundamental
solution introduced by Guenther, Thomann \cite{GT} for  the system (\ref{1.10}).

\smallskip
We will use the following technical lemmas:

\begin{lem}\label{repeating}
Let $\delta >0.$ Assuming  $z\in B_\delta(x) ,$ we have \begin{equation}\label{evident}|z|\ge|x|/2,\ \ \mbox{for}\ |x|\geq 2\delta,\end{equation} \begin{equation}\label{staudelta} s_\tau (z)^{-1} \le \mathfrak C\, (1+|x-z|)\, s_\tau(x)^{-1}\le \mathfrak C(\delta)\, s_\tau(x)^{-1}.\end{equation}  \end{lem}
{\bf Proof:} For $\ |x|\geq 2\delta $ we have $|z|\ge|x|-|x-z|\ge|x|-\delta\ge|x|/2,  $ i.e. the relation (\ref{evident}) is satisfied.  For the proof of (\ref{staudelta}) see  \cite[Lemma 4.8]{DK2}.

\begin{lem}[\mbox{\cite[Corollary 3.1]{DKN3}}] \label{corollary3.10}
Let $ j,k \in \{1,\, 2,\, 3\} ,\; \alpha , \beta \in \mathbb{N} _0^3$
with $| \alpha + \beta |\le 2,\\ \;y,z \in \mathbb{R}^3 ,\; t \in (0, \infty ) .$
Then
\begin{eqnarray*} &&
| \partial _y ^{ \alpha } \partial _z ^{ \beta } \Gamma _{jk} (y,z,t)|
\le
\mbox{$\mathfrak C $} \, (|y- \tau \, t \, e_1- e^{-t \cdot \Omega } \cdot z|^2+t)
^{-3/2-| \alpha + \beta |/2}.
\end{eqnarray*}
\end{lem}
\smallskip

\begin{lem}[\mbox{\cite[Theorem 2.19]{DKN2}}]\label{lemma2.60}
Let $S_1,S \in (0, \infty ) $ with $S_1<S,\; \nu \in (1, \infty )$. Then
\begin{eqnarray} \label{2.60.10} && \hspace{-3em}
\int_{ 0}^{ \infty } (|y- \tau \, t \, e_1- e^{-t \cdot \Omega } \cdot z|^2 +t)
^{-\nu}\; dt
\le
\mbox{$\mathfrak C$} (S_1,S,\nu ) \, \bigl(\, |y| \cdot s_{ \tau }(y) \,\bigr) ^{-\nu + 1/2}
\end{eqnarray}
for $ y \in B^c_{ S}, \; z \in \overline{ B_{S_1}}.$
\end{lem}
\smallskip

\begin{lem}[\mbox{\cite[Lemma 3.2]{DKN3}}] \label{lemma3.2}
Let $ j,k \in \{1,\, 2,\, 3\} .$ For $ \alpha , \beta \in \mathbb{N} _0^3 $ with
$| \alpha + \beta |\le 2,\; y,z \in \mathbb{R}^3 $ with $y\ne z,$ the function
$
(0, \infty ) \ni t \mapsto \partial _y ^{ \alpha }  \partial _z ^{ \beta  } \Gamma _{jk}
(y,z,t)         \in  \mathbb{R}
$
is integrable, the derivative
$ \partial _y ^{ \alpha }  \partial _z ^{ \beta  } \mbox{$\mathfrak Z$}   _{jk} (y,z)$
exists, and
\begin{eqnarray} \label{3.2.1} &&
\partial _y ^{ \alpha }  \partial _z ^{ \beta  } \mbox{$\mathfrak Z$}   _{jk} (y,z)
=
\int_{ 0} ^{ \infty } \partial _y ^{ \alpha }  \partial _z ^{ \beta }
\Gamma   _{jk} (y,z,t)\; dt.
\end{eqnarray}
Moreover, for $ \alpha ,\, \beta $ as before, the derivative
$ \partial _y ^{ \alpha }  \partial _z ^{ \beta  } \mbox{$\mathfrak Z$}   _{jk} (y,z)$
is a continuous function of $y,z \in \mathbb{R}^3 $ with $y\ne z$.
\end{lem}

\begin{lem}
\label{theorem3.20}
Let $S_1, S \in (0, \infty ) $ with $S_1 < S$, $\alpha , \beta \in \mathbb{N} _0^3
$
with
$
| \alpha + \beta |\le 2,\; 1\le j,k\le 3.$ Then
\begin{eqnarray*}
|\partial _y ^{ \alpha }  \partial _z ^{ \beta  } \mbox{$\mathfrak Z$}   _{jk} (y,z)|
\le
\mbox{$\mathfrak C$} (S_1,S) \, \bigl(\, |y| \cdot s_{\tau } (y) \,\bigr) ^{-1-| \alpha +
\beta |/2}
\ \ \ \mbox{for}
\
y \in B_S^c,\; z \in \overline{ B_{S_1}},&& \\
|\partial _y ^{ \alpha }  \partial _z ^{ \beta  } \mbox{$\mathfrak Z$}   _{jk} (y,z)|
\le
\mbox{$\mathfrak C$} (S_1,S) \, \bigl(\, |z| \cdot s_{\tau } (z) \,\bigr) ^{-1-| \alpha +
\beta |/2}
\ \ \ \mbox{for}
\
z \in B_S^c,\; y \in \overline{ B_{S_1}}.&&
\end{eqnarray*}
\end{lem}
\smallskip

{\bf Proof:} Lemma \ref{lemma2.60} - \ref{lemma3.2}.

\smallskip
\begin{lem}[\mbox{\cite[Theorem 3.1]{DKN1}}] \label{lemma2.70}
Let $ k \in \{0,1\},\;
R \in (0, \infty ) ,\; y,z \in B_R$ with $y\ne z.$
Then
\begin{eqnarray} && \nonumber
\int_{ 0} ^{ \infty } \bigl(\, |y- \tau \,\,  t \,\,  e_1 - e^{-t \,\,  \Omega } \cdot z|^2
+t \,\bigr) ^{-3/2-k/2}\; dt
\le
\mbox{$\mathfrak C$} (R) \,\,  |y-z|^{-1-k}
.
\end{eqnarray}

Due to Lemma \ref{corollary3.10}, this means for $y,z$ as above,
and for $ j,k \in \{1,\, 2,\, 3\} ,\;
\alpha \in \mathbb{N} _0^3$ with $| \alpha |\le 1$ that
\begin{equation*}
\left|\partial _y^{\alpha } \mbox{$\mathfrak Z$}(y,z)\right| + \left|\partial _z^{\alpha } \mbox{$\mathfrak Z$}(y,z)\right|\le
\mbox{$\mathfrak C$}(R) \,\,  |y-z|^{-1-| \alpha |}.
\end{equation*}

\smallskip
\end{lem}

\begin{lem}[\mbox{\cite[Lemma 4.1]{DKN3}}]
 \label{lemma4.1}
\ Let $j,k\in \{1,\, 2,\, 3\} ,\; g \in L^1( \partial \mbox{$\mathfrak D$} ),$
and put
\begin{eqnarray*} &&
F(y):= \int_{ \partial \mbox{$\mathfrak D$} }   \mbox{$\mathfrak Z$} _{jk}
(y,z) \, g(z)\; do_z \mbox{\ \ for}\;\;
y \in \overline{ \mbox{$\mathfrak D$}  }^c.\ \ \ \ \ \ \ \ \ \ \
\end{eqnarray*}
Then $F \in C^1( \overline{ \mbox{$\mathfrak D$} }^c)\; $ and
\begin{eqnarray} \label{4.1.10} &&
\partial_m F(y)= \int_{ \partial \mbox{$\mathfrak D$} } \partial y_m
 \mbox{$\mathfrak Z$} _{jk}
(y,z) \, g(z)\; do_z\ \ \  \hbox{for
$
1\le m\le 3,\;
\;
y \in \overline{ \mbox{$\mathfrak D$}  }^c.
$
}
\end{eqnarray}

\end{lem}
\begin{lem}

(\mbox{\cite[Lemma 4.2]{DKN3}})
 \label{lemma4.2}\ Let $j,k,l \in \{1,\, 2,\, 3\} ,\; \;
R>0,\; g \in L^1( B_R ),$
and put
\begin{eqnarray*} && \hspace{-1em}
F(y):= \int_{ B_R }  \partial_{z_l}\mbox{$\mathfrak Z$} _{jk}
(y,z) \, g(z)\; dz
\quad
\mbox{for}\;\;
y \in \overline{ B_R}^c.
\end{eqnarray*}
Then
$F \in C^1( \overline{ B_R }^c)$ and
\begin{eqnarray*} &&
\partial_m F(y)= \int_{ B_R } \partial_{y_m}\partial_{z_l}
\mbox{$\mathfrak Z$} _{jk}
(y,z) \, g(z)\; dz
\ \ \  \hbox{for}\;
\
y \in \overline{ B_R}^c,\, \, 1\le m\le 3.
\end{eqnarray*}

\end{lem}

\smallskip

\begin{lem} \label{theorem3.2}
Let $\,\gamma \in (1/4, \infty )$. Then there is a constant $C(\gamma)>0 $ such that for all $x\in \mathbb R ^3:$
\begin{equation*}
\int_{ \mathbb{R}^3 } \bigl[\, (1+|x-y|) \, s(x-y) \,\bigr]^{-3/2} \,\left[(1+|y|)s(y)^{}
\, \right]^{-\gamma} dy
\le
C(\gamma) \, (1+|x|)^{-c} \, s(x)^{-d}\,\ln^k(2+|x|),
\end{equation*}
where

$ c:=\begin{cases}
\gamma-1/2 & \hbox{if }\ \ \gamma\in(1/4,\,2] \\
3/2 & \hbox{if }\ \ \gamma\in(2,\,+\infty) \\
\end{cases}
$
$\ \  d:=\begin{cases}
\gamma & \hbox{if }\ \ \gamma\in(1/4,\,3/2] \\
3/2 & \hbox{if }\ \ \gamma\in(3/2,\,+\infty) \\
\end{cases}
$
\ \ \
$
k:=\begin{cases}
0 & \hbox{if }\ \ \gamma\ne2 \\
1 & \hbox{if }\ \ \gamma=2. \\
\end{cases}
$
%
%
\end{lem}
{\bf Proof:} See the proof of \cite[Theorem 3.2]{KNP}.

\begin{lem}
\label{convolution}
There exist a constant $C>0$ such that  for all $ x\in\mathbb R^3:$
\begin{equation*}\displaystyle \int_{\mathbb R ^3}\left[\left(1+|x-y|\right) \,s(x-y)\right]^{-2}\left[\left(1+|y|\right)\,s(y)\right]^{-2}\hbox{\rm d}y  \leq C\left[\left(1+|x|\right)\,s(x)\right]^{-2}\ln(2+|x|)
\end{equation*}
\end{lem}

{\bf Proof.}

Let us denote  $r^*:=\min(1,r),\, r\in\mathbb R$,\, $\eta^\alpha_\beta(x):=\left(1+|x-y|\right)^\alpha s(x-y)^\beta$, $x,y\in\mathbb R^N,\; N \in \mathbb{N},$ $\; N\ge 2,\;$$\alpha,\beta\in \mathbb R.$ In  \cite{KNP} the  inequalities of the type $\eta ^{-a}_{-b} \ast\eta^{-c}_{-d}\leq c\,\eta^{-e}_{-f}$ are studied in $\mathbb R^N.$  If $ |x |\le 1, $ the conditions  $ a+b^*+c+d>N $ and $ a+b+c+d^*>N$ ensure that the convolution $\eta ^{-a}_{-b} \ast\eta^{-c}_{-d} $ is bounded by $ C(a,b,c,d)$, and thus by $ C(a,b,c,d)\,\eta^{-e}_{-f}$,  (\cite[p.  73]{KNP}). So, we may consider the case $|x|\ge 1$. In that case   the whole space $\mathbb R^N$ is divided into sixteen regions $\Omega _i,\ 0\leq i\leq 15, $
 and the optimal choice of $e_i, f_i$   in the inequality $$\eta ^{-a_{}}_{-b_{}} \ast\eta^{-c_{}}_{-d_{}}\leq C_{i}\,\eta^{-e_{i}}_{-f_{i}}\ \ \hbox{in} \ \ \Omega _i,\ \ 0\leq i\leq 15,  $$ for given $a,b,c,d $\ is stated  in \cite[Tab.\,1,\,2]{KNP}  included in this paper as an appendix. Using the expressions of $e_i $ and $f_i$, we have to find $$e=\min_{i=0,1, ...15}e_i,\ \ \ e+f=\min_{i=0,1, ...15}(e_i+f_i).$$


Let us mention that $\eta^\alpha_\beta(x) \leq 2^{\gamma-\alpha }\,\eta^\gamma_\delta(x), $ for $x\in \mathbb R^N$ if $\alpha\leq\gamma $ and $\alpha+\beta\leq \gamma +\delta,$ $ \alpha, \beta, \gamma, \delta\in \mathbb R.$
 Using expressions of $e_i$ from  \cite[Tab.\,1,\,2]{KNP}, where  we put $N=3,$ we define $e$ as the minimum of the following values: \
\begin{eqnarray*} &c+{1\over 2}\min(0,\,a+b^*-3),\ &a+\textstyle{1\over 2}\min(0,\,c+d^*-3),\, \\
&a+c-2+{1\over 2}\min(0,1+b^*-a,1+d-c), \ &a+c-2+\textstyle{1\over 2}\min(0,1+b-a,1+d^{*}-c),
\\
&a+c+d-2+{1\over 2}\min(0,1-c-d),\ &a+c+b-2+\textstyle{1\over 2}\min(0,1-a-b), \\
&a+b+c+d-3+\textstyle{1\over 2}\min(0,3-2b-c-d),\ &a+b+c+d-3+\textstyle{1\over 2}\min(0,3-a-b-2d),
\\&a+b^{*}+c+d-3,\ &a+b+c+d^*-3.
\end{eqnarray*}
Substituting $a=b=c=d=2,\,  $ we get $ e=2.$
Analogously, using expressions of $e_i, \, f_i$ from  \cite[Tab.\,1,\,2]{KNP} we define $e+f$ as the minimum of the following values \
\begin{eqnarray*} &c+d+\min(0,\,a+b^*-3),\ &a+b+\min(0,\,c+d^*-3),\, \\
&a+b^{*}+c+d-3, \ &a+b+c+d^*-3.
\end{eqnarray*}
So, we get $ e+f=4, $ hence $ f=2.$

The logarithmic factor
of the type $ \ln(2+|x|)$ appears on $\Omega_0 $ because $a+b^*=3,$
and on  $ \Omega_1$ because $c+d^*=3$, see \cite[Tab.\,1,\,2]{KNP}. Regions $\Omega_2,\Omega_3\, $ and $\Omega_4$ also contribute  logarithmic factors, they are covered by the logarithmic factor of the mentioned type $ \ln(2+|x|)$.
\hfill $\Box$

\medskip
Starting point of our considerations will be the following theorem about the integrability and pointwise decays  of the velocity and its gradient, where the velocity is a solution of the rotational Navier-Stokes equations:


\begin {theo}[\mbox{\cite[Theorem 1.1]{DKN5}}]
 \label{theorem1.1}
Let $ \tau \in (0, \infty ) ,\; \omega \in \mathbb{R}^3 \backslash\{{\mathbf 0}\} ,\; \mbox{$\mathfrak D$} \subset \mathbb{R}^3 $
open and bounded. Take $ \gamma,\, S_1 \in (0, \infty ) ,\; p_0 \in (1, \infty ),\; A \in (2, \infty ), \; B \in [0, 3/2],\;
f: \mathbb{R}^3 \mapsto \mathbb{R}^3 $ measurable with $ \overline{ \mbox{$\mathfrak D$} } \subset B_{S_1},\;
A+\min\{B,1\} >3,\; A+B\ge 7/2,\;
f|B_{S_1} \in L^{p_0}(B_{S_1})^3,$
\begin{eqnarray} \label{T1.1.1}
|f(y)|\le \gamma  \,|y|^{-A} \, s_\tau(y)^{-B}
\quad \mbox{for}\;\; y \in B_{S_1}^c.
\end{eqnarray}
Let
$u \in L^6( \overline{ \mbox{$\mathfrak D$}} ^c)^3\cap W^{1,1}_{loc}( \overline{ \mbox{$\mathfrak D$} }^c)^3
,\;  \pi \in L^2_{loc}( \overline{ \mbox{$\mathfrak D$}} ^c)$, $\nabla u\in L^2( \overline{ \mbox{$\mathfrak D$} }^c)^9,\;$                                                                                                                                                                                                                                                                                                                                                                                                                                                                                                                                                                                                                                                                                                                                                                                                                                                                                                                                                                                                                                                                                                                                                                                                                                                                                                                                                                                                                                                                                                                                                                                                                                                                                                                                                                                                                                                                                                                                                                                                                                                                        and
\begin{eqnarray}\nonumber &&
\displaystyle\int_{ \overline{ \mbox{$\mathfrak D$} }^c} \biggl[\,
\nabla  u \cdot \nabla \varphi  \; +\;  \bigl(\, \tau \, \partial _1 u
 + \tau \, (u \cdot \nabla )u
\\ && \hspace{3em}
- (\omega \times z) \cdot \nabla  u +
\omega \times u \,\bigr) \cdot \varphi  -  \pi \, \mbox{div } \varphi - f \cdot \varphi \,\biggr]\; dz = 0,\ \ \  \mbox{div} \, u =0\qquad
\label{weakproblem}\end{eqnarray}

for $ \varphi \in C ^{ \infty } _0( \overline{ \mbox{$\mathfrak D$}   }^c)^3.$
Let $S \in (S_1, \infty ).$ Then
\begin{eqnarray} \label{1.4b}
| \partial ^{\alpha }u(y)| \le
D
\,  \bigl(\, |y| \, s_\tau(y) \,\bigr) ^{-1-| \alpha |/2}
\quad \mbox{for}\;\; x \in B_S^c,\;\; \alpha \in \mathbb{N} _0^3\; ,
\, \ |\alpha |\le 1,
\end{eqnarray}
with the constant $D$
depending on $\tau,\rho,\gamma ,\, S_1,\, p_0,\, A,\, B,\, \|f|B_{S_1}\|_1,\, u, \, \pi,\, S,$ and
on an arbitrary but fixed number $S_0 \in (0,S_1) $ with $\overline{ \mbox{$\mathfrak D$}}  \subset B_{S_0}$.
\end{theo}

 Let
$p \in (1, \infty ),\; q \in (1,2),\; f \in L^p_{loc}( \mathbb{R}^3 )^3
$
with
$  f|B_S^c \in L^q(B_S^c)^3 $
for some $ S \in (0, \infty ) .$

For  $ y \in \mathbb{R}^3,$ $j \in \{1,\, 2,\, 3\} $,
we set
\begin{eqnarray*} &&
\mbox{$\mathfrak R$}_j (f)(y):= \int_{ \mathbb{R}^3 } \sum_{k = 1}^3 \mbox{$\mathfrak Z$} _{jk}
(y,z) \, f_k(z)\; dz.\ \ \ \ \ \ \
\end{eqnarray*}
 According to \cite[Lemma 3.1]{DKN2}, the integral
appearing in the definition of $\mbox{$\mathfrak R$} _j(f)$  is well defined at least for almost every $y \in \mathbb{R}^3 $. If $ f $ is a function on $  \overline{ \mbox{$\mathfrak D$}}^c,$ the function $f$ in the previous definition is to be replaced by the extension of $ f $  by zero to $ \mathbb R^3.$

In order to derive the leading terms of the velocity and its gradient we are going to use the representation formula of a solution of the rotational Navier-Stokes equation:

%
\medskip
\begin{theo} \label{theorem5.20}
Let $u \in W^{1,1}_{loc} ( \overline{ \mbox{$\mathfrak D$} }^c)^3\cap
L^6( \overline{ \mbox{$\mathfrak D$} }^c)^3$ with
$ \nabla u \in L^2( \overline{ \mbox{$\mathfrak D$} }^c)^9.$
Let $\; p \in (1, \infty ),\;
q \in (1,2),\; f: \overline{ \mbox{$\mathfrak D$} }^c \mapsto \mathbb{R}^3 $ a function
with $f| \mbox{$\mathfrak D$} _T \in L^p( \mbox{$\mathfrak D$} _T)^3$ for
$T \in (0, \infty ) $ with $ \overline{ \mbox{$\mathfrak D$} } \subset B_T,\;
f|B_S^c \in L^q( B_S^c )^3$ for some $S \in (0, \infty ) $ with $ \overline{ \mbox{$\mathfrak D$}}
\subset B_S$.
Further assume that $  u| \mbox{$\partial\mathfrak D$}\in W^{2-1/p,p}( \partial{ \mbox{$\mathfrak D$} })^3  $ and $\; \pi: \overline{ \mbox{$\mathfrak D$} }^c \mapsto \mathbb{R} $ is a function
with $\pi| \mbox{$\mathfrak D$} _T \in L^p( \mbox{$\mathfrak D$} _T)$ for
$T$ as above.

Suppose that the pair $ (u,\pi ) $ is a weak solution
of the Navier-Stokes system with
Oseen and rotational terms, and with right-hand side $f$ in the sense of (\ref{weakproblem}).
Then $u \in W^{2,\,\min\{p,3/2\}}( \mbox{$\mathfrak D$} _T)^3,\;  \pi \in W^{1,\,\min\{p,3/2\}}( \mbox{$\mathfrak D$} _T)$
for any $T \in (0, \infty ) $ with $\overline{ \mbox{$\mathfrak D$} } \subset B_T$,
\begin{eqnarray} \label{5.100} && \hspace{-3em}
u_j(y)= \mbox{$\mathfrak R$} _j \bigl(\, f - \tau \,\,  (u \cdot \nabla )u \,\bigr) (y)
+ \mbox{$\mathfrak B$} _j(u,\pi) (y)
\quad \mbox{for}\;
j \in \{1,\, 2,\, 3\}  ,\;  \mbox{a.e.} \ y \in \overline{ \mbox{$\mathfrak D$} }^c,
\end{eqnarray}
where $\mbox{$\mathfrak B$} _j(u,\pi)$ is defined by
\begin{eqnarray} \label{4.199}
 \mbox{$\mathfrak B$} _j&&\hspace{-6.pt}(u,\pi )(y)
\\  \nonumber&&
\ := \int_{ \partial \mbox{$\mathfrak D$} }
 \sum_{k = 1}^3  \Bigl[
\sum_{l = 1}^3 \Bigl(  \mbox{$\mathfrak Z$}  _{jk} (y,z)
\,\,
\bigl(\, -\partial _lu_k(z) + \delta _{kl} \,\,  \pi(z) + u_k(z) \,\,
(\tau \,\,  e_1 - \omega \times z)_l \,\bigr)
\\ && \nonumber \hspace{4em}
+
\partial z_l \mbox{$\mathfrak Z$}  _{jk} (y,z)\,\,  u_k(z) \Bigr)
\,\,  n ^{( \mbox{$\mathfrak D$} )}_l(z)
\; +\;
E_{4j}(y-z) \,\,   u_k(z) \,\,  n_k ^{( \mbox{$\mathfrak D$} )}(z)
\Bigr]\; do_z
\end{eqnarray}
for $y \in \overline{ \mbox{$\mathfrak D$} }^c$.
\end{theo}
{\bf Proof:} See \cite[Theorem 4.1]{DKN6}, and its proof, as well as \cite[Theorem 4.4]{DKN2}. \hfill $\Box$
 \smallskip

In comparison with the linear case we will need some additional lemma:
\begin{lem} Let $ \, \phi \in W^{1,1}_{loc}(U)$ for some open set $U \subset \mathbb{R}^3$  and $ A\in \mathbb R^{3\times 3}$$\hbox{ such that}\ $ $  A^{-1}=A^T. \ \hbox{Then:}$
\[
\, \ A\nabla _z\left( \phi(Az)\right)= \nabla\phi(Az)
\]
\end{lem}
{\bf Proof:} Indeed:
\[ \frac{\partial }{\partial z_l}\left(\phi(Az)\right)  =\sum  _{k=1}^3\partial _k\phi(Az)\frac{\partial (Az)_{k} }{\partial z_l}=\sum  _{k=1}^3\partial _k\phi(Az) A_{k\, l}=\sum_{k=1}^3 A_{l\, k}^T\,\partial _k\phi(Az)\] i.e.
\(
\nabla _z\left( \phi(Az)\right)=A^T\, \nabla\phi(Az),
\) which gives the mentioned formula.$\hfill\Box$

\medskip
\medskip
\begin{cor}\label{corollaryperpartes}
In the situation of Theorem \ref{theorem1.1}, we get for $z \in \overline{ B_{S_1}}^c$ that \[ \sum_{l=1}^{3}\left(u_l\,\partial _lu\right)(\hbox{e}^{t\Omega}z)=\sum_{l=1}^{3}{\partial _l}\left( u_lu_{}\right)_{}(\hbox{e}^{t\Omega} z)\]

\[
=\sum_{l=1}^{3}
\sum_{k=1}^{3}
\left(
    \hbox{e}^{t\Omega}
\right)_{l k}
{\,\frac{\partial }{\partial z_k}}
\left[
    \left( u_lu_{}
    \right)
    (
       \hbox{e}^{t\Omega} z
     )
\right]
=
\sum_{l=1}^{3}
\left(  e^{t\Omega}{\,\nabla _z}
\right)_l
\left[\left(
  u_lu
\right)
(\hbox{e}^{t\Omega} z)\right].
\]
\end{cor}


\begin{lem}\label{continuous_differentiability}
In the situation of Theorem \ref{theorem1.1}, we have
\begin{eqnarray} \label{*} &&
\int_{ \overline{ \mbox{$\mathfrak D$} }^c}| \partial _x^{\alpha } \mbox{$\mathfrak Z$}(x,y)\,
[(u \cdot \nabla )u](y)|\; dy < \infty
\quad \mbox{for}\;\;
x \in \overline{ B_{S_1}}^c,\; \alpha \in \mathbb{N} _0^3\;\;\mbox{with}\; | \alpha |\le 1.
\end{eqnarray}
Moreover the function
$
\mbox{$\mathfrak V$}(x):=  \int_{ \overline{ \mbox{$\mathfrak D$} }^c}\mbox{$\mathfrak Z$}(x,y)\,
[(u \cdot \nabla )u](y)\; dy\;\; (x \in  \overline{ B_{S_1}}^c)
$
belongs to $C^1( \overline{ B_{S_1}}^c)^3,$ with
\begin{eqnarray} \label{**} &&
\partial ^{\alpha }\mbox{$\mathfrak V$}(x)
=
\int_{ \overline{ \mbox{$\mathfrak D$} }^c} \partial _x^{\alpha } \mbox{$\mathfrak Z$}(x,y)\,
[(u \cdot \nabla )u(y)|\; dy
\quad \mbox{for}\;\;
x,\, \alpha\; \mbox{as in (\ref{*})}.
\end{eqnarray}
\end{lem}

{\bf Proof:}
Let $ U \subset  \mathbb{R}^3 $ be open and bounded, with $\overline{ U} \subset  \overline{ B_{S_1}}^c.$
It is enough to show that (\ref{*}) holds for $x \in U$, that $\mbox{$\mathfrak V$}|U \in
C^1(U)^3,$ and (\ref{**}) is valid for $x \in U$.

Due to our assumptions on $U$, we may choose $R,\, S \in (S_1, \infty )$ such that
$\overline{B_{S}}\cap \overline{ U}= \emptyset$ and $\overline{ U}\subset B_R$. In particular we have
$\,\mbox{dist}( B_{S},U) >0$ and $\mbox{dist}(U,B_R^c)>0. $ This observation and Lemma \ref{theorem3.20} imply that
$
| \partial ^{\alpha } _x\mbox{$\mathfrak Z$}(x,y)|\le C_0$ for $x \in U,\; y \in B_{S} \backslash \overline{
\mbox{$\mathfrak D$} },\; \alpha \in \mathbb{N} _0^3$ with $| \alpha |\le 1,$ where $C_0$ is independent
of $x$ and $y$. We further observe that $(u \cdot \nabla )u \in L^{3/2}( \overline{ \mbox{$\mathfrak D$} }^c)^3,$
hence $(u \cdot \nabla )u| B_{S}\backslash \overline{ \mbox{$\mathfrak D$}} \in L^1(
B_{S}\backslash \overline{ \mbox{$\mathfrak D$} })^3$.
Lemma \ref{theorem3.20} and (\ref{1.4b}) yield that
$
| \partial ^{\alpha }_x \mbox{$\mathfrak Z$}(x,y)\, [(u \cdot \nabla )u](y)|
\le
C_{1} \cdot |y|^{-7/2-| \alpha |/2}
$
for $x \in U,\; y \in B_{R}^c$ , with $C_1$ again being independent of $x$ and $y$.
In view of the last statement of Lemma \ref{lemma3.2}, we may thus conclude by Lebesgue's theorem
that the function
\begin{eqnarray*}
y \mapsto  \partial ^{\alpha }_x \mbox{$\mathfrak Z$}(x,y)\, [(u \cdot \nabla )u](y),\; y \in
A:=( B_{S}\backslash \overline{ \mbox{$\mathfrak D$} })\cup B_{R}^c,
\end{eqnarray*}
is integrable for $x \in U,\; \alpha \in \mathbb{N} _0^3$ with $| \alpha |\le 1,$
that the function
\begin{eqnarray*}
\mbox{$\mathfrak V$}^{(I)}(x):= \int_{ A} \mbox{$\mathfrak Z$}(x,y)\, [(u \cdot \nabla )u](y)\; dy,
\;\; x \in U,
\end{eqnarray*}
belongs to $C^1(U)^3,$ and that
$
\partial ^{\alpha } \mbox{$\mathfrak V$}^{(I)}(x)
=
\int_{ A}\partial ^{\alpha }_x \mbox{$\mathfrak Z$}(x,y)\, [(u \cdot \nabla )u](y)\;dy
$
for $x,\, \alpha $ as before.

Let $\varphi \in C ^{ \infty } _0( \mathbb{R}^3 )$ with $0\le \varphi \le 1,\; \varphi |B_{1/2}=0,\;
\varphi |B_1^c =1,$
and define
$\varphi _{\delta }(x):= \varphi ( \delta ^{-1} x)$ for $x \in \mathbb{R}^3 ,\; \delta >0$.
Then
$
\varphi _{\delta }\in C ^{ \infty } _0( \mathbb{R}^3 ),\; 0\le \varphi _{\delta }\le 1,\; \varphi_{\delta } |B_{\delta /2}=0,\;
\varphi_{\delta } |B_{\delta }^c =1$ and $| \nabla \varphi _{\delta }(x)|\le  C \, \delta ^{-1} $
for $x \in \mathbb{R}^3 ,\; \delta >0$.

Using Lemma \ref{lemma2.70} and Theorem \ref{theorem1.1}, we see there are constants $C_2,C_3$ with
\begin{eqnarray}\nonumber && \hspace{-3em}
| \partial ^{\alpha }_x \bigl(\,  \mbox{$\mathfrak Z$}(x,y)\, \varphi _{\delta }(x-y) \,\bigr) \, [(u \cdot \nabla )u](y)|
\\&&
\le
C_2 (|x-y|  ^{-2} +\delta ^{-1} |x-y| ^{-1} )\,\chi_{( \delta /2, \, \infty )}(|x-y|)\label{***}
\le
C_3 \delta ^{-2}
\end{eqnarray}
for
$
x \in U,\; y \in B_{R}\backslash B_{S},\; \alpha \in \mathbb{N} _0^3
$
with $| \alpha |\le 1$.
In addition, if $y \in B_{R}\backslash B_{S},$
the function
$
x \mapsto \mbox{$\mathfrak Z$}(x,y)\, \varphi _{\delta }(x-y)  \, [(u \cdot \nabla )u](y),\; x \in U,
$
is continuously differentiable, as follows from Lemma 2.6. Now we may conclude from Lebesgue's theorem
that the function
$
y \mapsto
\partial ^{\alpha }_x \bigl(\,  \mbox{$\mathfrak Z$}(x,y)\, \varphi _{\delta }(x-y) \,\bigr) \, [(u \cdot \nabla )u](y),\;
y \in B_{R}\backslash B_{S},
$
is integrable for any $\delta >0,\; x \in U,\; \alpha \in \mathbb{N} _0^3$ with $| \alpha |\le 1,$
the function
\begin{eqnarray*}
\mbox{$\mathfrak B$} _{\delta }(x):=
\int_{ B_{R}\backslash B_{S}}  \mbox{$\mathfrak Z$}(x,y)\, \varphi _{\delta }(x-y)
\, [(u \cdot \nabla )u](y)\; dy,\;\;x \in U,
\end{eqnarray*}
belongs to $C^1(U)^3$ for any $\delta >0$, and
\begin{eqnarray*}
\partial ^{\alpha }\mbox{$\mathfrak B$} _{\delta }(x)
=
\int_{B_{R}\backslash B_{S}}
\partial ^{\alpha }_x \bigl(\,  \mbox{$\mathfrak Z$}(x,y)\,
\varphi _{\delta }(x-y) \,\bigr) \, [(u \cdot \nabla )u](y)\; dy
\end{eqnarray*}
for $\delta ,\, x,\, \alpha $ as before. Proceeding as in (\ref{***}), we further obtain
\begin{eqnarray*} &&
\int_{B_{R}\backslash B_{S}}
|\partial ^{\alpha }_x \bigl(\,  \mbox{$\mathfrak Z$}(x,y)\,
\varphi _{\delta }(x-y) - \mbox{$\mathfrak Z$}(x,y) \,\bigr) \, [(u \cdot \nabla )u](y)|\; dy
\\ &&
\le
C_4 \int_{B_{R}\backslash B_{S}} \chi_{B_{\delta }}(x-y) \,(|x-y| ^{-2} + \delta ^{-1} |x-y| ^{-1} )\; dy
\le
C_5 \delta
\end{eqnarray*}
for $x \in U,\; \alpha \in \mathbb{N} _0^3$ with $| \alpha |\le 1$, with the  $C_4,C_5$ denoting constants
independent of $\delta $ and $x$. Therefore, by an argument involving uniform convergence of $\mbox{$\mathfrak B$}
_{\delta }$ and $\nabla \mbox{$\mathfrak B$} _{\delta }$ for $\delta\downarrow  0$, we may conclude that the function
\begin{eqnarray*}
\mbox{$\mathfrak V$}^{(II)}(x):=
\int_{B_{R}\backslash B_{S}}
\mbox{$\mathfrak Z$}(x,y)\, [(u \cdot \nabla )u](y)\; dy, \;\; x \in U,
\end{eqnarray*}
belongs to $C^1(U)^3$, and
\begin{eqnarray*}
\partial ^{\alpha }\mbox{$\mathfrak V$}^{(II)}(x)
=
\int_{B_{R}\backslash B_{S}}
\partial _x^{\alpha }\mbox{$\mathfrak Z$}(x,y)\, [(u \cdot \nabla )u](y)\; dy
\quad \mbox{for}\;\; x \in U,\;
\alpha \in \mathbb{N} _0^3\;\mbox{with}\;| \alpha |\le 1.
\end{eqnarray*}
Since $ \mbox{$\mathfrak V$}(x)= \mbox{$\mathfrak V$}^{(I)}(x)+\mbox{$\mathfrak V$}^{(II)}(x)$
for $x \in U,$ the proof of the lemma is complete.
\hfill $\Box $

\section{ Leading term of the velocity and of its gradient}
{The aim of this part is to find the leading term of the velocity and its gradient for the Navier-Stokes  problem with rotation. Let us recall that} {the quantities}  \( \tau  ,\; \omega \) and the set $\mathfrak D $ were fixed in Section 2. We study the case $f$ has a compact support in $ \overline{ \mbox{$\mathfrak D$}} ^c. $ The result we will prove in the work at hand may be stated as:


\begin{theo}
 \label{mainLemma} Let $S_1 \in (0, \infty ) $ with $\overline{ \mbox{$\mathfrak D$} } \subset B_{S_1},\; p \in (1, \infty ),\;
f \in L^p( \overline{ \mbox{$\mathfrak D$} }^c)^3$ with $supp(f) \subset B_{S_1},\;
u \in L^6( \overline{ \mbox{$\mathfrak D$} }^c)^3 \cap W^{1,1}_{loc}( \overline{ \mbox{$\mathfrak D$} }^c )^3$ with
$\nabla u \in L^2( \overline{ \mbox{$\mathfrak D$} }^c)^9 $ and $
u| \partial \mbox{$\mathfrak D$} \in W^{2-1/p,\,p}( \partial \mbox{$\mathfrak D$} )^3,\;
\pi \in L^2_{loc}( \overline{ \mbox{$\mathfrak D$} }^c)$ with $\pi| \mbox{$\mathfrak D$} _{S_1}
\in L^p( \mbox{$\mathfrak D$} _{S_1})$.
Suppose that the pair $ (u,\pi ) $ is a weak solution
of the Navier-Stokes system with
Oseen and rotational terms, and with right-hand side $f$ in the sense of (\ref{weakproblem}).
Then there are coefficients $ \beta  _1, \beta  _2, \beta  _3 \in \mathbb{R} $
and functions $ \mbox{$\mathfrak F$} _1, \mbox{$\mathfrak F$} _2, \mbox{$\mathfrak F$} _3
\in C^1( \overline{ B_{S_1}}^c)$ such that for $ j \in \{1,2, 3\} ,\;
\alpha \in \mathbb{N} _0^3$ with $ | \alpha |\le 1,\;
x \in \overline
{ B_{S_1}}^c$,
%
%
%
\begin{equation}\label{leading}\partial ^{ \alpha }u_j(x)=\left\{
\sum_{k = 1}^3 \beta  _k \, \partial ^{ \alpha }\mbox{$\mathfrak Z$}_{jk} (x,0)
+
 \left(\int_{ \partial \mbox{$\mathfrak D$} } u \cdot n ^{( \mbox{$\mathfrak D$} )}  \, do_z
\right) \, \partial ^{ \alpha }E_{4j}(x)\right\}
+
\partial ^{ \alpha }\mbox{$\mathfrak F$} _j(x),\end{equation}and if   $\, S \in (S_1, \infty ), $
$ x \in B_S^c$,
\begin{equation}\label{T4.2.20aNL}
| \partial ^{\alpha } \mbox{$\mathfrak F$} (x)|
\le
\mbox{$\mathfrak C$} \,(|x|\, s_{\tau }(x) )^{-3/2-| \alpha |/2} \ln(2+|x|),
\end{equation}
where \mbox{$\mathfrak C$} depends on  $\tau ,\, \omega ,\, p,\, S_1,\, S,$ certain norms of $u,\, \pi $ and
$f$, and on the constant $D$ from (\ref{1.4b}).

\end{theo}

In Theorem \ref{mainLemma}, the estimate presented in \cite[Theorem 3.14]{DKN6} for the linear case is extended
to the nonlinear one. Note that by \cite[(3.9)]{GT}, the function $ \mbox{$\mathfrak Z$}(x,0)$
in the leading term on the right-hand side of (\ref{leading})
corresponds to the time integral of a fundamental solution of the evolutionary Oseen system multiplied by a rotation depending on time.

{\bf Proof of Theorem \ref{mainLemma} }The term  of  (\ref{leading}) contained in braces $\{\dots \}  $ we will call "the leading term",  term $ \mathfrak F$ we will call "the remainder". From Theorem
\ref{theorem5.20}  we have \begin{eqnarray} \label{5.100NL}
\hspace{-3em} u_j(x)= \mbox{$\mathfrak R$} _j \bigl(\, f - \tau
\,(u \cdot \nabla )u \,\bigr) (x) + \mbox{$\mathfrak B$} _j(u,\pi)
(x), \ j \in \{1,\, 2,\, 3\},\;   \mbox{for}\; \ \mbox{a.e.} \; x
\in \overline{ \mbox{ $\mathfrak D $} }^c, \ \ \ \ \end{eqnarray}

where $ \mbox{$\mathfrak B $} _j(u,\pi)$ was defined in
(\ref{4.199}).
{


We put
\begin{eqnarray*} &&
\beta _k:=\beta^{(I)}_k-\tau\beta ^{(II)}_k \\ &&\beta ^{(I)}_k:=\int_{ B_{S_1}} f_k(y)\; dy
\\ && \hspace{2em}
+
\int_{ \partial \mbox{$\mathfrak D$} } \sum_{l = 1}^3 \bigl(\,
- \partial _lu_k(y) + \delta _{kl} \, \pi(y)
+
u_k(y) \, ( \tau \, e_1- \omega \times y)_l \,\bigr)
\, n^{ ( \mbox{$\mathfrak D$}  )}_l(y)
\;do_y
\\ && \beta ^{(II)}_k:=\int_{\partial \overline{\mathfrak D}}  \,   \sum_{m=1}^3\, (n_{m}u_{m}\,u_k) (y)\, \hbox{d}o_{y}
\end{eqnarray*}
for $1\le k\le 3.$ By the definition of $\beta _k$ the leading term  in formula $(\ref{leading})$ is determined. Because  $(\ref{leading})$   is in fact rearrangement of formula ($\ref{5.100NL}  $), we now define  the value $\mathfrak F_j $ as the difference of the right-hand side of the representation formula $ (\ref{5.100NL})$\ minus the leading term. We will distinguish  $\mathfrak F^{(I)} $ coming from the linear terms and  $\mathfrak F^{(II)}$  arising from the non-linear part, i.e. from $ \,\mbox{$\mathfrak R$} _j \bigl(\, \,(u \cdot \nabla )u \,\bigr):$
\begin{eqnarray}&&
\hspace{-0.5em}\mbox{$\mathfrak F$} _j(x)
:= \mbox{$\mathfrak F^{(I)}$} _j(x)
-\tau\,\mathfrak F^{(II)} _j(x), \nonumber\\ &&
\hspace{-0.5em}\mbox{$\mathfrak F^{(I)}$} _j(x):=\int_{ B_{S_1}}
\Bigl(
\sum_{k = 1}^3 \bigl[\, \bigl(\,  \mbox{$\mathfrak Z$}_{jk}
(x,y) - \mbox{$\mathfrak Z$} _{jk} (x,0) \,\bigr)
\, f_k(y) \,\bigr]
\Bigr)
\; dy
 \label{T4.2.45}
\nonumber\\ && \hspace{0.0em} \nonumber
+
\int_{ \partial \mbox{$\mathfrak D$} }
\sum_{k = 1}^3 \Bigl(
\bigl(\,  \mbox{$\mathfrak Z$}_{jk}(x,y) - \mbox{$\mathfrak Z$}_{jk} (x,0) \,\bigr)
\\ && \hspace{7.5em} \nonumber
\cdot
\sum_{l = 1}^3 \bigl(\, - \partial _lu_k(y) + \delta _{kl} \, \pi(y)
+
u_k(y) \, ( \tau \, e_1 - \omega \times y)_l \,\bigr)
\,
n^{ ( \mbox{$\mathfrak D$}  )}_l(y)
\\ && \hspace{0.0em} \nonumber
+ \bigl(\, E_{4j}(x- y) - E_{4j}(x) \,\bigr)
\, u_k(y) \, n ^{( \mbox{$\mathfrak D$} )}_k(y)  \Bigr)  \; do_y
+
\int_{ \partial  \mbox{$\mathfrak D$} }
\sum_{k,l = 1}^3
\partial y_l\, \mbox{$\mathfrak Z$}_{jk}
(x,y) \, u_k(y) \, n ^{( \mbox{$\mathfrak D$} )}_l(y)\; do_y\,,
\end{eqnarray}
\begin{eqnarray}\label{ZminusZ}
\hspace{-0.5em}\mbox{$\mathfrak F^{(II)}$} _j(x):=   \int_{\overline {\mathfrak D}^c}\sum_{k,l = 1}^3 \mbox{$\mathfrak Z$}_{jk}(x,y)\left (u_{l}\,\partial_l u_k\right) (y)\hbox{d}y
-\int_ {\partial\mbox{$\mathfrak D $}}\sum_{k,l=1}^3 \mbox{$\mathfrak Z$}_{jk}(x,0)(n_{l}u_{l}\,u_k) (y)\, \hbox{d}o_{y}
\quad
\end{eqnarray}
for $x \in \overline{ B_{S_1}}^c,\; 1\le j\le 3.$ Then by (\ref{5.100NL}) we get (\ref{leading}).

\smallskip


The assertion of the theorem  will be proved in four steps:\medskip

 {\bf 1.\ Estimates and continuity  of $\partial^\alpha \mathfrak F^{(I)},$ where $|\alpha|=0$ or $|\alpha| =1 $.
}

By exactly the same proof as given in [5, p. 473-474] for [5, Theorem 1.1], we obtain
that  $ \mbox{$\mathfrak F$}^{(I)} \in C^1( \overline{ B_{S_1}}^c)^3 $ and
\begin{eqnarray*}
| \partial ^{\alpha } \mbox{$\mathfrak F$}^{(I)}(x)|
\le
\mbox{$\mathfrak C$} (\|f\|_1 + \| \nabla u| \partial \mbox{$\mathfrak D$} \|_1
+  \| \pi| \partial \mbox{$\mathfrak D$} \|_1 +  \|u| \partial \mbox{$\mathfrak D$} \|_1 )
(|x|\, s_{\tau }(x))^{-3/2-| \alpha |/2}
\end{eqnarray*}
for $S\in(S_1,\infty), $ $x \in \overline{ B_S}^c,\, \alpha \in \mathbb{N} _0^3$ with $| \alpha |\le 1$,
with \mbox{$\mathfrak C$} depending on  $\tau ,\, \omega ,\, p,\, S_1$ and $S.$

\medskip {\bf 2.  $C^1$-continuity of $\mathfrak F ^{II}$. }

By Lemma \ref{continuous_differentiability}, the function $\mbox{$\mathfrak F$} ^{(II)}\in C^1( \overline{ B_{S_1}}^c)^3$, and first-order derivatives may be moved into the volume integral
appearing on the right-hand side of (\ref{ZminusZ}).

\medskip
{ \bf 3.\ Estimates  of $\partial^\alpha \mathfrak F^{(II)}$: first steps.
}

Let $x \in B_S^c$. Recalling that
\begin{eqnarray}
    \mathfrak F^{(II)}(x)=\int_{\overline {\mathfrak D}^c}\sum_{l = 1}^3 \mbox{$\mathfrak Z$}(x,y)\left( u_{l} \, \partial _l u\right) (y)\hbox{d}y\,\ -\int_ {\partial\mbox{$\mathfrak D $}}\sum_{l=1}^3\, \mbox{$\mathfrak Z$}(x,0)(n_{l}u_{l}\,u) (y)\, \hbox{d}o_{y},
\end{eqnarray}

we apply firstly the integration by parts and and then split the resulting volume integral in an integral $\mbox{$\mathfrak B$}_R\,$   on the bounded domain $ B_R\setminus \overline{\mathfrak D}$ and integral $\mathcal E_R$ on the exterior  domain $ (B_R)^c,$ where $R=(S_1+S)/2$. Thus $ \mathfrak F^{II}(x) $ becomes
\[\int_{\partial {\mathfrak D}}  \,   \sum_{l=1}^3\, \left[\mathfrak Z_{} (x,y)-\mbox{$\mathfrak Z$}(x,0)\right](n_{l}u_{l}\,u) (y)\, \hbox{d}o_{y}-\left\{\int_{B_R\setminus \overline{\mathfrak D}}+\int_{(B_R)^c}\right\}  \sum_{l=1}^3\partial _{y_l}\,  \mathfrak Z_{} (x,y)(u_{l}\,u) (y)\, \hbox{d}y \]
\begin{equation}={\mathcal S}_{\partial\mathfrak D} -{\mathcal B_R}-\mathcal E_R.\label{split100}\end{equation}

Of course, here and in similar situations in the following, a partial integration has to be performed
first on a bounded domain, where $B_T \backslash (\overline{ B_R}\cup B_{\epsilon }(x)) $
with $T>\max\{2R,\, 2|x|\},\; 0< \epsilon $ is a good choice for such a domain.
In the next step we let $\epsilon $ tend to zero. This passage to the limit may be handled
by referring to Lemma \ref{lemma2.70} and (\ref{1.4b}).
Finally we let $T$ tend to infinity.
The surface integral on $\partial B_T$ which came up in the partial integration
then vanishes, as follows from Lemma \ref{theorem3.20} and (\ref{1.4b}). The same references imply that all
the volume integrals involved tend to integrals on $B_R^c$ when $T\to \infty $.

In  volume integral $\mathcal E_R $ over the exterior domain $(B_R)^c $ we use   firstly the definition of $\mathfrak Z $, (\ref{3.2.1}) and the Fubini's     theorem, and then the domain invariant transformation $y= e^{t\Omega} z$ for fixed $t>0.$ The reason why we use the mentioned transformation is that we would like to avoid a
periodic  term in the right-hand side of  (\ref{trick}):  \[\mathcal E_R=\int_{_{(B_R)^c}}  \sum_{l=1}^3\partial _{y_l}\,  \mathfrak Z_{} (x,y)(u_{l}\,u) (y)\, \hbox{d}y= \int_{0}^{\infty}\int_{{ (B_R)}^c}   \sum_{l=1}^3\partial _{y_l}\,   \Gamma_{} (x,y,t) [u_{l}\,u] (y)\, \hbox{d}y\, \hbox{d}t \]    \[= \int_{0}^{\infty}\int_{{ (B_R)}^c}   \sum_{l=1}^3 \partial _{y_l}\,   \Gamma_{} (x,y,t)_{|y={e}^{t\Omega}  z}[u_{l}\,u] (e^{t\Omega}z)\, \hbox{d}z\, \hbox{d}t \]    Finally we split the the exterior domain of integration $B^c_R$ on two domains: Let $\delta$ be a sufficiently small positive number comparing to $1, R$ and $S-S_1$,    f.e. $\delta:=\mbox{min}\{1,(S-S_1)/2,R/2\}=\mbox{min}\{1,S-R,R/2\}.$
Note that $\overline{B_\delta(x)}\subset B_R^c.$ We obtain: \[{\mathcal E_R}= \left\{\int_{0}^{\infty}\int_{  B_\delta{(x)}}+\int_{0}^{\infty}\int_{ {(B_R)}^c\setminus B_\delta{(x)}} \right\}    \sum_{l=1}^3\partial _{y_l}\, \Gamma_{} (x,y,t)_{|y={e}^{t\Omega}  z} [(u_{l}\,u) (\hbox{e}^{t\Omega}z)\, \hbox{d}z\, \hbox{d}t\]\[ =\mathcal V_\delta+\mathcal V_{R,\delta} \]

 Substituting the expression of $\mathcal E_R $ into ($\ref{split100}$) we get:
\begin{equation} {\mathfrak F^{II}}={\mathcal S}_{\partial\mathfrak D} -{\mathcal B_R-\mathcal V_\delta-\mathcal V_{R,\delta}}\label{1symbolFII}\end{equation}


$\partial _x^{\alpha}\mathcal B_R, \ \partial _x^{\alpha} \mathcal  S_{\partial\mathfrak D}:$ \  Estimating the behavior of the first two terms and their derivatives $ \partial^\alpha_x $ for $|\alpha|=0,1 $, we get the following estimate:
\begin{equation} |\partial _x^{\alpha}\mathcal B_R|+|\partial _x^{\alpha} \mathcal  S_{\partial\mathfrak D}|\le \mathfrak C(S_{1},S)\, (|x|s_\tau(x))^{-3/2-|\alpha|/2}, \ x\in B^c_S\label{help1}\end{equation}Indeed, from  Lemma \ref{theorem3.20} for $y \in B_R,\, x\in B^c_S$:\ \begin{equation}\label{number1}|\partial _x^{\alpha}\partial_ {y_l}\,  \mathfrak Z_{} (x,y)|\le \mathfrak C(S_1,S)\,(|x|s_\tau(x))^{-3/2-|\alpha|/2},  \end{equation}
\begin{equation}\label{number2}  \left|\partial _x^{\alpha}\big(\mathfrak Z_{} (x,y)- \mathfrak Z_{} (x,0)\big)\right| =\left|\sum _{k=1}^3 \,\partial _x^{\alpha}\partial_ {y_k} \,\mathfrak Z_{} (x,\theta y) \,y_k\right|\le \mathfrak C(S_1,S)\, (|x|s_\tau(x))^{-3/2-|\alpha|/2}\end{equation} for some   $0\leq\theta \leq 1.$
 So, with Lemma \ref{lemma4.2} and (\ref{number1}) \[ |\partial _x^{\alpha}{\mathcal B}_R| \le \left|\int_{B_R\setminus \overline{\mathfrak D}}  \sum_{l=1}^3 \partial _x^{\alpha}\partial _{y_l}\,  \mathfrak Z_{} (x,y)(u_{j}\,u) (y)\, \hbox{d}y \right| \le \] \[\le \mathfrak C(S_1,S)(|x|s_\tau(x))^{-3/2-|\alpha|/2} \int_{B_R\setminus \overline{\mathfrak D}}  \sum_{l=1}^3\big|(u_{l}\,u) (y)\big|\hbox{d}y\le  \mathfrak C(S_{1},S)\, (|x|s_\tau(x))^{-3/2-|\alpha|/2}, \,  \] because $|u|^2$  is $L^1$-integrable on bounded domain $B_R\setminus \overline{\mathfrak D}.$
Similarly,   we have with (\ref{4.1.10}) and  (\ref{number2}):
\begin{eqnarray*}  |\partial _x^{\alpha} \mathcal  S_{\partial\mathfrak D}|
%
%
&\le& \mathfrak C(S_1,R)\,(|x|s_\tau(x))^{-3/2-|\alpha|/2} \int_{\partial {\mathfrak D}}    \sum_{l=1}^3\big|(n_{l}u_{l}\,u) (y)\,\big| \hbox{d}o_{y}
%
\\
%
&\le&  \mathfrak C(S_{1},R)\, (|x|s_\tau(x))^{-3/2-|\alpha|/2} .
\end{eqnarray*}

\medskip
{ \bf 4.\ Estimates  of $\partial^\alpha \mathfrak F^{(II)}$ for $\alpha=0 $.
}

 $\mathcal V_\delta$:\ \ For  the estimation of this term  we
use Lemma \ref{corollary3.10}
for the first order derivatives of $\Gamma:
$ We have (for $x\not=\hbox{e}^{t\Omega }z    $) \begin{equation} \left|\partial_{y_j}\Gamma_{} (x,y,t)_{|y={e}^{t\Omega}  z}\right|\le \mathfrak C \left( \left| x-\tau t \hbox{e}_{1} -z\right|^{2}+t \right)^{-2}.
\label{trick}\end{equation}
 From  Theorem \ref{theorem1.1}
 we have for $y\in B_R^c $ \begin{equation*}\label{u_squared}|u(y)|^{2}  \le \mathfrak C(R)\left(|y|\,s_\tau(y)\right)^{-2}.
 \end{equation*}
If $z\in(B_R)^c  $ then $\, \hbox{e}^{t\Omega}  z\in(B_R)^c,$ we get:   \begin{equation}\label{usquared} |u(\hbox{e}^{t\Omega}  z)| ^{2}\le \mathfrak C(R)\left(|\hbox{e}^{t\Omega}  z|\,s_\tau(\hbox{e}^{t\Omega}  z)\right)^{-2}=\mathfrak C(R)\left(| z|\,s_\tau( z)\right)^{-2} \end{equation}\ Since $\overline{ B_{\delta }(x)} \subset B_R^c$, we thus get due to (\ref{evident}), (\ref{staudelta})
\begin{eqnarray}
\label{3.13a}
|u(e^{t \Omega }z)| \le \mbox{$\mathfrak C$} (R, \delta )\, (|x| \, s_{\tau }(x)) ^{-2}  \quad \mbox{for}\;\;
z \in \overline{ B_{\delta }(x)}.
\end{eqnarray}  
So, we have
\[ \left|\mathcal V_\delta\right|=\left|\int_{0}^{\infty}\int_{ B_\delta(x)}   \sum_{j=1}^{3}\partial _{y_j}\,   \Gamma (x,y,t)_{|y=e^{t \Omega  }z}\left[\left( u_ju_{}\right)_{}(e^{t\Omega} z)\right] \mbox{d}{z}\, \hbox{d}t\right|\] \[\le \mathfrak C(R) \int_{ B_\delta(x)}  \int_{0}^{\infty}  \left( \left| x-\tau t \hbox{e}_{1} -z\right|^{2}+t \right)^{-2}\left(| x|\,s_\tau( x)\right)^{-2}\,\hbox{d}t \, \hbox{d}{z}\,  \]
\[ \le \mathfrak C(R, \delta)\left(| x|s_\tau( x)\right)^{-2}_{} \int_{ B_\delta(x)}  |x-z|^{-2}\, \hbox{d}{z}\le \mathfrak C(R,\delta)\left(| x|s_\tau( x)\right)^{-2} ,\, \]  where the integral with respect to variable $ t $ is estimated using Lemma \ref{lemma2.70},
choosing in its application
$y:=x-z,\, $ $z:=0.$

 $ \mathcal V_{R,\delta}:$  Similarly as in the previous case using (\ref{trick}) and (\ref{usquared}):





\[|\mathcal V_{R,\delta}|\le  \left|\int_{0}^{\infty}\int_{ B_R ^c\setminus B_\delta(x)}   \sum_{j=1}^{3}\partial _{y_j}\,   \Gamma (x,y,t)_{|y=e^{t \Omega  }z}\left[\left( u_ju_{}\right)_{}(e^{t\Omega} z)\right] \mbox{d}{z}\, \hbox{d}t\right| \]

\[ \leq \ \mathfrak C(R) \int_{0}^{\infty}\int_{B_R ^c\setminus B_\delta(x)}   \left(\left |x-\tau t\mbox{e}_1-z\right|^2+t\right)^{-2}\left(\left| z\right|\, s_\tau(z)\right)^{-2}\hbox{d}z\, \hbox{d}t\]

Now, the integral with respect to $ t$      can be estimated using Lemma \ref{lemma2.60},
 $y:=x-z,$ $\  z:=0:$\[ \int _0 ^{\infty}\left(\left|x-\tau t\mbox{e}_1-z\right|^2+t\right)^{-2}\hbox{d}t\le \mathfrak C (S_{1},S)\left(|x-z|s_\tau (x-z)\right)^{-3/2},\ z\in B_R ^c\setminus B_\delta(x) \]

\begin{equation*}|\mathcal V_{R,\delta}|
%
\le \mathfrak C (S_{1},S)\int _{B_R ^c\setminus B_\delta(x)} \left(|x-z|s_\tau (x-z)\right)^{-3/2}\left(| z|\, |s_\tau( z)|\right)^{-2}\hbox{d}z
  \label{convolution1}\end{equation*}

\[ \le \mathfrak C (S_{1},S)\ln(2+|x|)\left(| x||s_\tau( x)|\right)^{-3/2}.\hbox{}\] The last inequality
follows from Lemma \ref{theorem3.2} ($\gamma=2$).



\medskip
\bf 5.\ Estimates  of $\partial^\alpha \mathfrak F^{(II)}$ for $ |\alpha| =1 $.

\rm
Let us mention that  $S, \,S_1, \, R, \,\delta $ are the same as in the previous section, so $\overline{B_\delta(x)}\subset B^c_R. $  The aim of this part is to find the leading term of the gradient of velocity for the Navier-Stokes  problem with rotation:
The difference with the previous case is that we cannot apply the integration by parts over the whole domain $\overline{\mathfrak D}^c$ because we have to protect the neighbourhood $B_{\delta}(x) $ due to singularities of the second order derivatives of $\mathfrak Z.$ On the other hand, to avoid some technical difficulties,  we are able to handle the integrals with respect to $ t $ only in domains invariant with respect to the transformation   $y=\mbox{e}^{t\Omega}z, \ t>0.$ These facts causes  some additional computations.
So,  we use Lemma \ref{continuous_differentiability}, split the domain of integration into the bounded part $ B_R\setminus \overline{{\mathfrak D}}^c$ and the exterior domain $(B_R)^c, $ and we apply the integration by parts firstly only on the bounded domain:      \begin{eqnarray}\displaystyle
&&\partial^\alpha _x\mbox{$\mathfrak F^{(II)}_j$} (x)=  \int_{\overline {\mathfrak D}^c}\sum_{k,l = 1}^3 \partial^\alpha _x\mbox{$\mathfrak Z$}_{jk}(x,y)\left (u_{l}\,\partial_l u_k\right) (y)\hbox{d}y-\int_ {\partial\mbox{$\mathfrak D $}}\sum_{k,l=1}^3\,\partial^\alpha _x \mbox{$\mathfrak Z$}_{jk}(x,0)(n_{l}u_{l}\,u_k) (y)\, \hbox{d}o_{y}\nonumber
\\ &&
=\left\{\int_{B_R\setminus \overline{\mathfrak D}}+\int_{(B_R)^c}\right\}\sum_{k,l = 1}^3 \partial^\alpha _x\mbox{$\mathfrak Z$}_{jk}(x,y)\left (u_{l}\,\partial_l u_k\right) (y)\hbox{d}y \nonumber\\ &&\qquad\qquad\qquad-\int_ {\partial\mbox{$\mathfrak D $}}\sum_{k,l=1}^3\,\partial^\alpha _x \mbox{$\mathfrak Z$}_{jk}(x,0)(n_{l}u_{l}\,u_k) (y)\, \hbox{d}o_{y}  \nonumber\end{eqnarray}
\begin{eqnarray}\ \ \ =\int_{\partial \overline{\mathfrak D}}   \,   \sum_{k,l=1}^3\left[\partial^\alpha_x \mathfrak Z_{j k} (x,y)-\partial^\alpha_x \mathfrak Z_{j k} (x,0)\right](n_{l}u_{l}\,u_{k}) (y)\, \hbox{d}o_{y}\qquad\qquad\qquad\qquad\qquad\qquad\qquad\label{up} \nonumber\\   \nonumber-\int_{B_R\setminus \overline{\mathfrak D}}  \sum_{k,l=1}^3\partial {y_l}\,  \partial ^\alpha _x \mathfrak Z_{jk} (x,y)(u_{l}\,u_k) (y)\, \hbox{d}y+\int_{\partial B_R}  \sum_{k,l=1}^3\,    \partial ^\alpha _x \mathfrak Z_{jk} (x,y)(u_{l}\,u_k) (y)\, \frac{y_{l}} {R}\hbox{d}o_{y}\\ +\int_{(B_R)^c}\sum_{k,l = 1}^3 \partial^\alpha _x\mbox{$\mathfrak Z$}_{jk}(x,y)\left (u_{l}\,\partial_l u_k\right) (y)\hbox{d}y\qquad\qquad\nonumber  \qquad\qquad\qquad\qquad\qquad\qquad\qquad
 \end{eqnarray}So, we get:
\begin{equation}\partial^\alpha _x\mbox{$\mathfrak F^{(II)}$} (x)=\partial^\alpha\mathcal S_{\partial\mathfrak D} -\partial^\alpha{\mathcal B_R}+\mathcal S'_{R}+\mathcal E'_R\qquad\qquad \label{symbolFII} \end{equation}
 Evaluation of  the last term in (\ref{up}) with (\ref{3.2.1}):

 \[\mathcal E'_R(x)_j=\int_{(B_R)^c}\sum_{k,l = 1}^3 \partial^\alpha _x\mbox{$\mathfrak Z$}_{jk}(x,y)\left (u_{l}\,\partial_l u_k\right) (y)\hbox{d}y\qquad\qquad\qquad\qquad\qquad\]\[
=\int_{(B_R)^c}\int _0^{+\infty}\sum_{k,l = 1}^3 \partial^\alpha _x\mbox{$\Gamma$}_{jk}(x,y,t)\left (u_{l}\,\partial_l u_k\right) (y)\hbox{d}y\,\hbox{d}t\qquad\qquad \]The domain of integration of $\mathcal E'_R$  is $(B_R)^c.$ This exterior domain is invariant with respect to the transformation $y={e}^{t\Omega}z, \ t>0.$ We use the same transformation to avoid periodic terms as in the case $|\alpha|=0:$
\[{\mathcal E}'_R(x)_j=\int _0^{+\infty}\int_{(B_R)^c}\sum_{k,l = 1}^3 \partial^\alpha _x\mbox{$\Gamma$}_{jk}(x,\hbox{e}^{\tau\Omega}z,t)\left (u_{l}\,\partial_l u_k\right) (\hbox{e}^{\tau\Omega}z)\hbox{d}z\,\hbox{d}t\ \ \ \ \ \ \]
Unlike the case $|\alpha|=0,$ the mentioned transformation is used {\it before} the  integration by parts.
 We split the domain of integration into two domains $B_\delta(x)$ and $(B_R)^c \setminus B_\delta(x). $ In the integral over  the unbounded  domain    we apply the identity from Corollary \ref{corollaryperpartes} and integrate by parts:
\[ \mathcal E'_R(x)_j=\int _0^{+\infty}\int_{ B_\delta(x)} \sum_{k,l = 1}^3 \partial^\alpha _x\mbox{$\Gamma$}_{jk}(x,{e}^{\tau\Omega}z,t)\left (u_{l}\,\partial_l u_k\right) ({e}^{\tau\Omega}z)\hbox{d}z\,\hbox{d}t\] \[
+\int_{0}^{\infty}
\int_{ \partial B_\delta(x)}
\sum_{k, l=1}^{3}
   \partial^\alpha _x  \Gamma_{j k } (x,{e}^{t\Omega}  z,t)
\left[
   \left(
      u_lu_{k}
   \right)
 ({e}^{t\Omega} z)
\right]
\left(
   {e}^{t\Omega}{\,(x-z)/ {\delta}}
\right) _l
\hbox{d} o_{z}\, \hbox{d}t
\]
\[
+\int_{0}^{\infty}
\int_{ \partial B_R}
\sum_{k, l=1}^{3}
  \partial^\alpha _x\Gamma_{j k} (x,{e}^{t\Omega}  z,t)
 \left[
    \left(
       u_lu_{k}
    \right)_{}
  ({e}^{t\Omega} z)
\right]
\left(
   {e}^{t\Omega}{\,( -z)/R}
\right) _l \hbox{d}o_{z}\, \hbox{d}t\] \[-\int_{0}^{\infty}\int_{{ B_R^c\setminus B_\delta(x)}}   \sum_{k,l=1}^{3}
\left(
   \hbox{e}^{t\Omega}{\,\nabla _z}
\right)_l
\partial^\alpha _x  \Gamma_{j  k}
(x,{e}^{t\Omega}  z,t)
\left[
   \left(
     u_lu_{k}
   \right)
   ({e}^{t\Omega} z)
\right] \mbox{d}z\, \hbox{d}t\]

\[=\, \big({\mathcal U}_{\delta}\big)_j+\big({\mathcal S}_{\delta}'\big)_j +\big(- \mathcal S '_{R}\big)_j+  \big({\mathcal U}_{R,\delta }\big)_j\ \ \ \ \ \ \ \ \ \ \ \ \ \ \ \ \]
\noindent
Substituting the expression of $\mathcal E'_R(x)$  into  ($\ref{symbolFII}$) and using (\ref{3.2.1}), we get finally:  \begin{equation} \partial^\alpha _x\mbox{$\mathfrak F^{(II)}$} (x)=\partial^\alpha\mathcal S_{\partial\mathfrak D} -\partial^\alpha{\mathcal B_{ R}}+\mathcal U_{\delta}+{\mathcal S}_{\delta}' + \cal U_{R,\delta} \label{finalFII} \end{equation}
Now we will estimate all terms of (\ref{finalFII}) for  $|\alpha|=1$:

$\partial^\alpha\mathcal S_{\partial\mathfrak D}, \ \partial^\alpha{\mathcal B_{ B_R}}$: \  From  (\ref{help1}) we know that \( |\partial _x^{\alpha} \mathcal  S_{\partial\mathfrak D}|+|\partial _x^{\alpha}\mathcal B_R|\le \mathfrak C _1(S_{1},S)\, (|x|\,s_\tau(x))^{-2}.\)

\smallskip
$ {\mathcal U}_{\delta}$: Estimates of this term are completely analogous to the evaluation of $ {\mathcal V}_{\delta }$  in the  case $|\alpha|=0$. Only difference is that  from  Theorem \ref{theorem1.1}:
$|u(y)| |\nabla u(y)|\le \mathfrak C(S)\left(|y|\,s_\tau(y)\right)^{-5/2}$ \, for \, $y\in (B_R)^c $:  We get $$ |{\mathcal U}_{\delta}|\le  \mathfrak C(S_{1},S)\left(|x||s_\tau(x)|\right)^{-5/2}.$$

${\mathcal S}_{\delta}'$: From Lemma \ref{corollary3.10}
for the first order derivatives of $\Gamma$,  $ x\not=\hbox{e}^{t\Omega }z    $\,:

\begin{equation*}
    \left|
        \partial^\alpha _x  \Gamma (x,\hbox{e}^{t\Omega}  z,t)
    \right|
        \le \mathfrak C
    \left(
    \left|
         x-\tau t \hbox{e}_{1} -z
    \right|^{2}+t \right)^{-2}.
\end{equation*}

 By (\ref{3.13a})
 \[\,|u(\hbox{e}^{t\Omega}  z)| ^{2}\le  \mathfrak C(R,\delta)\, \left(| x|\,s_\tau( x)\right)^{-2}\ \,  \hbox{for} \ \,
     z\in \partial B_\delta(x).
\]

It is also clear that
$ \left| e^{t \Omega} \,(x-z)/ {\delta} \right|=1 $  for
$ z\in B_\delta (x).$

So, we have
\[\left|{\mathcal S}_{\delta}'\right|
\le \mathfrak C(R)
\int_{ \partial B_\delta (x)}  \int_0 ^\infty
\left(
  \left|
       x-\tau t \hbox{e}_{1} -z
  \right|^2+t
\right)^{-2}
\left(
  |x|\, s_\tau(x)
\right)^{-2}\,\hbox{d} t \, \hbox{d} o_z
\]

\[ \le \mathfrak C_{}(R, \delta)
   \left(
     | x|\,s_\tau( x)
   \right)^{-2} \int_{\partial B_\delta(x)}  |x-z|^{-2}\,\hbox{d}o_{z}\le \mathfrak C_{}(R,\delta)\left(| x|\,s_\tau( x)\right)^{-2} \, \]

   where the integral with respect to variable $ t $ is estimated using Lemma \ref{lemma2.70},
$y:=x-z,\, $ $z:=0. $
 So, the integral ${\mathcal S}_{\delta}'$ belongs to the remainder.

  $ {\cal U}_{R,\delta}:$ We shall use Lemma \ref{corollary3.10},
  for the evaluation of the second order derivatives of the function $ \Gamma$:
\[
\left|
   \left(
     \hbox{e}^{t\Omega}{\,\nabla _z}
   \right)_j
 \partial^\alpha_ x  \Gamma_{} (x,\hbox{e}^{t\Omega}  z,t)
\right|
\le
\mathfrak C
\left(
   \left|
      x-\tau t\mbox{e}_1-z
   \right|^2+t
\right)^{-5/2} \] The integral with respect to $ t$ of the  right-hand side can be estimated using Lemma \ref{lemma2.60}
choosing $y, z$ from  the lemma
by the following way: $y:=x-z, \ z:=~0.   $ \, Hence: \[ \int _0 ^{\infty}\left(\left|y-\tau t\mbox{e}_1-z\right|^2+t\right)^{-5/2}\hbox{d}t\le \mathfrak C (S_{1},S)\,\left(|x-z|\, s_\tau (x-z)\right)^{-2}\]
Using (\ref{usquared}), we find
\[
|{\cal U}_{R,\delta}| \le
\left| \
\int _0 ^{\infty}
\int_{{ B_R^c\setminus B_\delta(x)}}
\sum_{j=1}^{3}
\left(
    \hbox{e}^{t\Omega}{\,\nabla _z}
\right)_j
  \partial x_m  \Gamma_{} (x,\hbox{e}^{t\Omega}  z,t)
\left[
   \left(
     u_ju_{}
   \right)_{}
  (\hbox{e}^{t\Omega} z)
\right] \mbox{d}z\, \hbox{d}t
\right|
\]

\begin{equation*}
\le \mathfrak C (S_{1},S)
\int _{B_R^c\setminus B_\delta(x)}
\left(
  |x-z|s_\tau (x-z)
\right)^{-2}
\left(
   | z||s_\tau( z)|
\right)^{-2}\hbox{d}z  \label{convolution12}\end{equation*}
\[ \le \mathfrak C (S_{1,}S)\left(| x|\,s_\tau( x)\right)^{-2}\mbox{ln}(2+|x|).\] The last inequality
we get by Lemma \ref{convolution}.
\hfill $\Box $

\medskip
Remark: So, finally we get that the leading term of $ \partial ^\alpha u_j\,$ is expressed as in the linear case
\[
\sum_{k = 1}^3 \beta  _k \, \partial _x^{ \alpha
}\mbox{$\mathfrak Z$}_{jk} (x,0),
\]
where $\beta =\left(\beta _1,
\beta _2, \beta _3 \right) $ contains additionally the term
$
\int_{\partial {\mathfrak D}}
\sum_{j=1}^3(n_{j}u_{j}\,u) (y)\, \hbox{d} o_{y}.
$

\bigskip
\begin{cor}
{Let $U\subset{\mathbb R}^3$ be open and bounded, $S_1\in (0,\infty)$ with $\overline{U}\subset B_{S_1}$, $p\in(1,\infty),$  $f\in L^p(U)^3$ with $ {supp}(f)\subset B_{S_1}$. Let $u\in W^{1,1}_{loc} (\overline{U}^c)^3$ with $u\in L^6(\overline{U}^c)^3,$  $\nabla u\in L^2(\overline{U}^c)^9$,  $\pi\in L^2_{loc}(\overline{U}^c)$, and suppose that   $u,\pi$ satisfy (\ref{weakproblem}) 
(weak form of rotational Navier-Stokes system, as in
Theorem \ref{theorem1.1}) with $ U $ in the place of $\mathfrak D$.
Let $S_0\in  (0, S_1) $ with $\overline{U}\subset B_{S_0}$. Then the conclusions of Theorem \ref{mainLemma} hold, with $\mathfrak D$
replaced by $B_{S_0}$ .
}
\end{cor}
{\bf Proof}: Obviously $(u\cdot\nabla)u\in L^{3/2}(\overline{U}^c)^3 $ so that $ f -(u\cdot\nabla)u \in L_{loc}^{{\min}\{p, 3/2\}}(\overline{U}^c)^3$.
Therefore, by interior regularity of the Stokes system, as stated in \cite[Theorem IV.4.1]{Ga1},  we have $u\in W^{2,\,\min\{p ,\, 3/2\} } _{loc} (\overline{U}
^c)^3 $ and $\pi\in W^{1,\,\min\{p ,\, 3/2\} } _{loc} (\overline{U}
^c)$;   also see the
proof of \cite[Theorem 5.5]{DKN2}. Now the corollary follows from Theorem \ref{mainLemma}, with $p_0 =
\min(p, 3/2)$ and $B_{S_0}$ in the place of $\mathfrak{D}$.
\hfill $\Box $




\medskip
{\bf {Acknowledgments}}

 The research of \v S.~N.  was supported
by Grant Agency of Czech Republic P201-13-00522S  and by RVO
67985840.}

\newpage

\tabcolsep1mm

\setbox0=\vbox{\noindent\hsize 26cm
\thispagestyle{empty}
\hspace*{-25.5cm}\vbox{
\scriptsize
\begin{tabular}[t]{|l|l|l|l|l|l|l|l|l|l|l|c|}
\multicolumn{12}{c}{} \\[-1.5cm]
\multicolumn{12}{c}{\normalsize Tab.1, see \cite{KNP} \ \ $N=3$} \\[5mm]\\
\hline
 &  &  &  &  & &  &   & & & &    \\
 &   &  &  &  &  & &   & & & &    \\
\mm Dom.\mm\mm & $t$ & $\varrho$ & $ r $ & $\tilde t $ & $ \tilde \varrho $
& $\tilde r$
& $ \eta_{-b}^{-a}({\bf y} ) $ & $ \eta_{-d}^{-c} ( {\bf x}-{\bf y} ) $
& $ e$ & $ f $ &  log.   \\
 &  &  &  &  & &  &   & & & &  factors
\Linn
 &  &  &  &  & &  &   & & & &    \\
$\Omega_0$  &  &  &  & $\sim x_1 $  & $\sim |{\bf x}'|$ & $\sim |{\bf x}| $
&   & & &  &  $ \ln R (b=1 \wedge a<2) \vee $   \\
  &  &  &  &  & &  &   & & & &  $\vee (b \neq 1 \wedge a+b^*=3)$  \\
 & $\mm (-\frac{1}{8}R^\upsilon ; \frac{1}{8}R^\upsilon) \mm$
& $\mm (0;\frac{1}{8}R^\upsilon)\mm$
& $\mm (0;\frac{1}{8}R^\upsilon) \mm$ & $R$
& $ R^\upsilon$ & $R$
& $ \mm r^{-a}(1+s ({\bf y}))^{-b} \mm $
& $ \mm R^{-c-2\sigma d} \mm $
&
 $ \mm c+\frac{1}{2}\min(0,a+b^*-3) \mm $
& $\mm d+\frac{1}{2}\min(0,a+b^*-3) \mm $
& $ \mm \ln^2 R(b=1 \wedge a=2) \mm $  \\
 &  &  &  &  & &  &   & & & &    \\
\hline
& & & & & & & & & & & \\
 &$\sim x_1 $  & $\sim |{\bf x}'|$ & $\sim |{\bf x}| $ & & & & & & & &
   $ \mm \ln R (d=1 \wedge c<2) \vee \mm $  \\
 $\Omega_1$ &  &  &  & &  &   & & & &
 & $\mm  \vee (d\neq 1 \wedge c+d^*=3)  \mm $   \\
 & $R$ & $ R^\upsilon$ & $R$
 & $\mm (-\frac{1}{8}R^\upsilon ; \frac{1}{8}R^\upsilon) \mm$
&  $\mm(0;\frac{1}{8}R^\upsilon)\mm$
&  $\mm(0;\frac{1}{8}R^\upsilon)\mm$
&  $R^{-a-2\sigma b}$
&  $\mm \tilde r^{-c}(1 + s ( {\bf x - y} ))^{-d}\mm $
&  $\mm a+\frac{1}{2}\min(0,c+d^*-3)\mm $
&  $\mm b+\frac{1}{2}\min(0,c+d^*-3)\mm$
&  $\mm  \ln^2 R \; (c=2 \wedge d=1) \mm $       \\
& & & & & & & & & & &  \\
\hline
  &  &  &  &  & & & & & & & $\mm \ln \frac{R}{1+s}(\min(1+b^*-a, \mm $
\\
 & $\sim r$  &   & $\sim t$ &  $\mm (-\frac{1}{8}R^\upsilon ;\mm $
& $\sim |{\bf x}'|$ &  & $\mm r^{-a} \quad \varrho < \sqrt t \mm $
& & $\mm a+c-2+\frac{1}{2}\min\mm $
& $\mm b^*+d-1-\frac{1}{2}\min \mm $
& $\mm  1+d-c)=0\wedge b \neq1) \mm $    \\
$\Omega_2$ &  &  &  &  & &  &   & & &
& $\mm \ln_{+} s \cdot \ln \frac{R}{1+s}\mm $     \\
 &  &  &  &  & &  &   & & & & $\mm  ( b=1 \wedge 1+d=0) \mm $  \\
 & $\mm (R^\upsilon ; R) \mm $
& $\mm (0;\frac{1}{8}R^\upsilon) \mm$
&  $\mm  (R^\upsilon ; R) \mm $
&  $\mm  R-R^\upsilon) \mm $
& $R^\upsilon $
& $\mm R +R^\upsilon-r \mm $
&  $\mm  r^{-a+b} \varrho^{-2b} \quad \varrho > \sqrt t \mm $
& $ \mm \tilde r^{-c+d} R^{-2d \upsilon}\mm $
& $\mm (0,1 + b^*-a,1+d-c) \mm $
&  $\mm (0,1 + b^*-a,1+d-c) \mm $
& $ \mm (\ln_{+}s \quad a<2)(\ln R \quad a>2)\mm $\\
 &  &  &  &  & &  &   & & & &
    $\mm \begin{array}{ll}
         (\ln R\; \ln\frac{R}{1+s} & a=2) \wedge b=1
         \end{array}
       \mm$ \\
 \hline
& & & & & & & & & & & \\
 &  &  &  &  & &  &   & & &
 & $\mm \begin{array}{ll}
         &\ln\frac{R}{1+s}(\min(1+b-a, \mm
        \end{array}
          $\\
 &$\mm (-\frac{1}{8}R^\upsilon; \mm $
&$\sim |{\bf x}'|$
&  &$\sim \tilde r$  & &$\sim \tilde t$  &
& $ \mm \tilde r^{-c} \quad \tilde \varrho< \sqrt \tau \mm $
& $\mm a+c-2+\frac{1}{2}\min \mm $
&  $\mm b+d^*-1-\frac{1}{2}\min\mm $
&  $\mm 1+d^*-c)=0   \wedge d \neq 1) \mm $  \\
 & & & & & & & & & &
& $\mm  \ln_{+}s \cdot \ln \frac{R}{1+s} (d=1 \wedge \mm $ \\
$\Omega_3$
&$\mm R-R^\upsilon) \mm $
&$R^\upsilon$
& $\mm R+R^\upsilon- \tilde r \mm$
& $\mm (R^\upsilon ; R) \mm $
& $\mm (0;\frac{1}{8}R^\upsilon) \mm $
& $\mm (R^\upsilon ; R) \mm $
& $\mm r^{-a+b}R^{-2b\upsilon} \mm $
& $ \mm \tilde r^{-c+d} \tilde \varrho^{-2d} \quad \tilde \varrho
> \sqrt \tau \mm $
& $\mm (0,1+b-a,1+d^*-c) \mm $
& $\mm (0,1+b-a,1+d^*-c) \mm $
&  $\mm  \wedge 1+b-a=0) \mm $  \\
 &  &  &  &  & &  &   & & & &  $\mm (\ln_{+}s \; c<2)(\ln R \; c>2) \mm $ \\
 &  &  &  &  & &  &   & & & & $\mm (\ln R \frac{R}{1+s} \ c=2)
 \wedge d=1  \mm  $  \\
\hline
& & & & & & & & & & & \\
$\Omega_4$ &  &   &  & &  & &  &
& see $ \Omega_2, \Omega_3 \mm $
& see $ \Omega_2, \Omega_3 \mm $
& see $ \Omega_2, \Omega_3 \mm $     \\
& & & & & & & & & & & \\
\hline
 & & & & & & & & & & & \\
 &  &  &  &  & &  & & & & &   $ (\ln_{+}s \quad b=1) \cdot$   \\
 $\Omega_5$ & $\sim r$  &  & $\sim t$  & $\sim  - \tilde r $ &
 & $\sim | \tilde t|$
 & $\mm  R^{-a} \quad \varrho < \sqrt t \mm $
 & $\mm  |\tilde t |^{-c-d}\mm $
 & $ \mm a+c+d-2+\frac{1}{2} \mm $
 & $ \mm b^{*}-1-\frac{1}{2}\min \mm $
 & $\mm \cdot (\ln \frac{R}{1+s}\ c+d=1) \mm $     \\
 & $R$ &  $\mm (0;R^\upsilon) \mm$
 & $R$  & $ \mm (-R;-R^\upsilon) \mm $
 &  $ \mm (0;R^\upsilon)\mm$
 &  $ \mm (R^\upsilon;R) \mm $
 &  $\mm R^{-a+b} \varrho^{-2b} \quad \varrho> \sqrt t\mm$
 & &  $ \mm (0,1-c-d) \mm $
 &   $\mm (0,1-c-d) \mm $  &    \\
\hline
& & & & & & & & & & & \\
 &  &  &  &  & &  & &
 &
&
 &   $( \ln_{+}s \quad d=1) \cdot$   \\
$\Omega_6$ &$  \sim -r$  &  &$\sim|t|$  & $\sim \tilde r$
&  & $\sim \tilde t$   & $|t|^{-a-b}$
&  $ \mm R^{-c} \quad \tilde \varrho < \sqrt {\tilde t} \mm$
&  \mm  $a+b+c-2+\frac{1}{2}\min \mm$
&  $ \mm d^{*}-1-\frac{1}{2} \min \mm $
&   $\mm \cdot (\ln \frac{R}{1+s} \quad a+b=1) \mm $   \\
 & $\mm (-R;-R^\upsilon) \mm$
 &  $\mm (0;R^\upsilon) \mm$
 &  $ \mm(R^\upsilon;R)\mm$
 & $R$
  &  $\mm (0;R^\upsilon) \mm$
  &$R$  &
& $ \mm R^{-c+d} \tilde \varrho^{-2d} \quad \tilde \varrho > \sqrt
{\tilde t} \mm$
& $\mm (0,1-a-b) \mm $&$\mm (0,1-a-b) \mm $   & \\
\hline
 &  &  &  &  & &  & & & & &    \\
$\Omega_7$
&  & $\sim \tilde \varrho$
&  &
&{  $\sim \varrho, \tilde r$ }
&$\sim \tilde \varrho$
& & &  $\mm a+b+c+d-3+ \mm$
&  $\mm -\frac{1}{2} \min(0,3-2b- \mm$
&  $\mm \ln  \frac{R}{1+s} \quad 2b+c+d=3 \mm$   \\
  &$R$
&  $\mm (| \tilde t|;R) \mm$
& $R$
&  $\mm (-R;-R^\upsilon) \mm$
&  $ \mm (| \tilde t|;R)\mm$
&  $(\mm | \tilde t|;R) \mm$
&  $\mm R^{-a+b}  \varrho^{-2b}\mm$
& $\mm\varrho^{-c-d} \mm$
& $\mm \frac{1}{2}\min(0,3-2b-c-d) \mm$
& {  $ -c-d)$ }
&    \\
\hline
 &  &  &  &  & &  & & & & &    \\
$\Omega_8 $
&   &   $\mm\sim \tilde \varrho,r \mm$
& $\sim \varrho$
&   & $\sim \varrho$
&  & &
&  $ \mm a+b+c+d-3+ \mm $
&  $ \mm -\frac{1}{2}\min(0,3-a-b-2d) \mm$
&  $\mm  \ln  \frac{R}{1+s} \quad a+b+2d=3 \mm$   \\
 & $\mm (-R;-R^\upsilon) \mm$
&  $ \mm (| t|; R) \mm$
&  $ \mm(|t|; R) \mm$
&  $R$
&   $\mm (|  t|; R) \mm$
& $R$
&  $ \mm \tilde \varrho^{-a-b} \mm$
&  $\mm R^{-c+d} \tilde \varrho^{-2d} \mm$
&  $ \mm \frac{1}{2}\min(0,3-a-b-2d) \mm$
&
&  \\
\hline
\end{tabular}
}
}

\box0

\newpage
\tabcolsep1mm

\setbox1=\vbox{\noindent\hsize 26cm
\thispagestyle{empty}
\hspace*{-25.62cm}\vbox{
\scriptsize
 \begin{tabular}[t]{|l|l|l|l|l|l|l|l|l|l|l|c|}
\multicolumn{12}{c}{} \\[1cm]
\multicolumn{12}{c}{\normalsize Tab.2, see \cite{KNP} \ \ $N=3$} \\[5mm]
\hline
 &  &  &  &  & &  &   & & & &    \\
 &  &  &  &  & &  &   & & & &    \\
\mm Dom.\mm & $t$ &  $\varrho$ &  $ r $ &  $\tilde t $ &  $ \tilde \varrho $
&  $\tilde r$
&  $ \eta_{-b}^{-a}({\bf y} ) $ & $ \eta_{-d}^{-c} ( {\bf x}-{\bf y} ) $
& $ e$ & $ f $ &  log. factors     \\
 &  &  &  &  & &  &   & & & &
\Linn
 &  &  &  &  & &  &   & & & &    \\

$\Omega_9$               
&  $ \begin{array}{l}
  \mm\sim r \\           
 \mm R \mm
  \end{array} $
&  $\begin{array}{l}
  \mm \sim \tilde \varrho \\         
   \mm (R^\upsilon ; |\tilde t|)
   \end{array} $\mm
& $ \begin{array}{l}
  \mm  \sim \tilde t \\                   
    \mm R
    \end{array} $  \mm
&  $ \begin{array}{l}
  \mm  \sim -\tilde r \\               
   \mm (-R;-R^\upsilon)
   \end{array}$\mm
&  $\begin{array}{l}
  \mm \sim  \varrho \\                
   \mm (R^\upsilon ;|\tilde t|)
   \end{array}$\mm
& $ \begin{array}{l}
   \mm \sim |\tilde t| \\             
  \mm (R^\upsilon;R)
   \end{array}$ \mm

&  $ R^{-a+b} \varrho^{-2b}$ \mm         
&  $|\tilde t|^{-c-d}$ \mm              
& \mm $ \begin{array}{lll}
     b>1 & \mbox{see} & \Omega_5 \\   
     b<1 & \mbox{see} & \Omega_7 \\
     b=1 & \mbox{see} & \Omega_5
   \end{array}$ \mm
&                                    
&  $ \begin{array}{l}
  \mm  \left ( \ln \frac{R}{1+s}\;c+d<1) \right .
    \\  \\                           
  \mm  \left . (\ln^2\frac{R}{1+s}\; c+d=1) \right ) \wedge b=1 \\ \\
    \end{array} $\mm
\\
\hline


$\Omega_{10}$               
&$ \begin{array}{l}
 \mm \sim - r \\          
  \mm (-R;-R^\upsilon)
  \end{array} $\mm
& $\begin{array}{l}
 \mm  \sim \tilde \varrho \\         
 \mm  (R^\upsilon ; |t|)
   \end{array} $   \mm
& $\begin{array}{l}
    \sim |t| \\                 
 \mm   (R^\upsilon ; R)
    \end{array}  $\mm
& $ \begin{array}{l}
 \mm \sim \tilde r \\          
 \mm R
  \end{array} $\mm

& $\begin{array}{l}
 \mm  \sim \varrho \\           
 \mm  (R^\upsilon ; |t|)
    \end{array} $\mm
&  $\begin{array}{l}
 \mm  \sim \tilde t \\         
 \mm  R
    \end{array}$\mm
&$ |t|^{-a-b}$   \mm            
& $R^{-c+d} \; \tilde \varrho^{-2d}$\mm   
& $ \begin{array}{lll}
 \mm    d>1 & \mbox{see} & \Omega_6 \\   
 \mm    d<1 & \mbox{see} & \Omega_8 \\
 \mm    d=1 & \mbox{see} & \Omega_6
   \end{array}$  \mm
&                               
& $ \begin{array}{l}
 \mm   \left (\ln \frac{R}{1+s}\;a+b<1) \right .
    \\  \\                           
 \mm   \left . (\ln^2\frac{R}{1+s}\; a+b=1) \right ) \wedge d=1 \\ \\
    \end{array} $\mm
\\
\hline

 &  &  &  &  & &  &   & & & &    \\
$\Omega_{11}$
& \mm$(R;\infty)$\mm
& \mm$ (0;\infty)$\mm
&$\begin{array}{l}
 \mm   \sim \tilde r \\
 \mm   (R; \infty)
    \end{array}$\mm                
&\mm$ (-\infty ;-R^\upsilon)$\mm
&$ (0;\infty)$\mm
& $\begin{array}{l}
 \mm   \sim  r \\
 \mm   (R; \infty)
    \end{array}$   \mm             
& $(1+r)^{-a} (1+s({\bf y}))^{-b}$ \mm
&$ r^{-c-d}$\mm
&\mm$a+b^*+c+d-3>0$\mm
&$0$
&\mm$ \ln R\quad b=1$\mm
\\
\hline
 &  &  &  &  & &  &   & & & &    \\
 &  &  &  &  & &  &   & & & &    \\
$\Omega_{12}$
&\mm$(-\infty ; -R^\upsilon)$\mm
&\mm$(0; \infty)$\mm
&$\begin{array}{l}
 \mm   \sim \tilde r\\
 \mm   (R; \infty)
    \end{array}$ \mm               
& \mm $   (R; \infty)$\mm
&$   (0; \infty)$\mm
& $\begin{array}{l}
 \mm   \sim  r\\
 \mm   (R; \infty)
    \end{array}$ \mm               
& $\tilde r^{-a-b}$\mm
&\mm $(1+\tilde r)^{-c} (1+s({\bf x-y}))^{-d}$
&\mm$ a+b+c=d^*-3>0$\mm
&$0$
&\mm$ \ln R \quad d=1$\mm
\\
 &  &  &  &  & &  &   & & & &    \\
\hline
 &  &  &  &  & &  &   & & & &    \\
$\Omega_{13}$
&
& $\begin{array}{l}
 \mm \sim r, \tilde r \\
 \mm (R; \infty)
    \end{array}$  \mm              
& $\begin{array}{l}
 \mm \sim \tilde r \\
 \mm (R; \infty)
    \end{array}$  \mm              
&
& $\begin{array}{l}
 \mm \sim \varrho \\
  \mm(R; \infty)
    \end{array}$ \mm               
&$\begin{array}{l}
 \mm \sim  r \\
 \mm (R; \infty)
    \end{array}$ \mm               
& &
&\mm see $\Omega_{11}$ and $\Omega_{12}$\mm
&\mm see $\Omega_{11}$ and $\Omega_{12}$   \mm
&\mm see $\Omega_{11}$ and $\Omega_{12}$   \mm
\\
 &  &  &  &  & &  &   & & & &    \\
\hline
 &  &  &  &  & &  &   & & & &    \\
$\Omega_{14}$
& $\begin{array}{l}            
 \mm  \sim r \\
 \mm  R
    \end{array}$  \mm              
&$\begin{array}{l}
 \mm  \sim \tilde \varrho \\
 \mm  (R^\upsilon ; R)
    \end{array}$   \mm
&$\begin{array}{l}
 \mm   \sim t \\                    
 \mm   R
    \end{array}$\mm
&\mm$ (-\frac 18R^\upsilon ; \frac{R}{2} )$
& $\begin{array}{l}
 \mm   \sim\varrho  \\
 \mm   (R^\upsilon;R)
    \end{array}$  \mm              
&\mm$ \sim |\tilde t |+\tilde \varrho$\mm
&$ R^{-a+b} \; \varrho^{-2b}$\mm
&$\begin{array}{l}
    (\tilde t + \tilde \varrho)^{-c+d} \;      
   \tilde \varrho^{-2d} \; \tilde t>0\\  \\
    \tilde \varrho^{-c-d} \; \tilde t<0
    \end{array}$\mm
& $\begin{array}{ll}               
 \mm   \mbox{see} & \Omega_2, \Omega_3\\
               & \Omega_7\\
               & \Omega_2\\
              & a+b+c+d-3
    \end{array}$\mm
&  $\begin{array}{l}       
 \mm   b+d>1\wedge 1+d-c>0\\
 \mm   1+d-c<0\\
 \mm   1+d-c<0\wedge b+d>1\\
 \mm   0\quad \mbox{otherwise}
    \end{array}$\mm
&$\begin{array}{l}
 \mm  \ln\frac R{1+s}\; b+d=1\wedge 1+d-c>0\\
      \\                                 
 \mm  \ln^2\frac R{1+s} \; b+d=1 \wedge 1+d-c=0
    \end{array}$\mm
\\
 &  &  &  &  & &  &   & & & &    \\
\hline
 &  &  &  &  & &  &   & & & &    \\
$\Omega_{15}$
& \mm $ (-R^\upsilon ;\frac R2)$\mm
&$\begin{array}{l}
 \mm   \sim \tilde \varrho\\           
 \mm   (R^\upsilon;R)
    \end{array}$\mm
& $\sim |t|+\varrho$\mm
& $\begin{array}{l}
  \mm \sim r\\
  \mm R
    \end{array}$\mm               
& $\begin{array}{l}
  \mm \sim \varrho\\
  \mm (R^\upsilon ;R)
    \end{array}$\mm                
& $\begin{array}{l}
  \mm  \sim t\\
  \mm  R
    \end{array}$  \mm              
& $\begin{array}{ll}
   (t+\varrho)^{-a+b}\;\varrho^{-2b} &t>0\\
   \varrho^{-a-b} & t<0
    \end{array}$  \mm                   

& $R^{-c+d}\; \tilde\varrho^{-2d}$\mm

& $\begin{array}{ll}                 
  \mm  \mbox{see} & \Omega_2, \Omega_3\\
               & \Omega_8\\
               & \Omega_3\\
              & a+b+c+d-3
    \end{array}$\mm
&  $\begin{array}{l}               
  \mm  b+d>1\wedge 1+b-a>0\\
  \mm  1+b-a<0\\
  \mm  1+b-a=0\wedge b+d>1\\
  \mm  0\quad \mbox{otherwise}
    \end{array}$\mm
& $\begin{array}{l}
 \mm  \ln\frac R{1+s}\; b+d=1\wedge 1+b-a>0\\    
                                           \\
 \mm  \ln^2\frac R{1+s} \; b+d=1 \wedge 1+b-a=0
    \end{array}$
\\
 &  &  &  &  & &  &   & & & &    \\
\hline
\end{tabular}
}
}
\box1


\medskip
\begin{thebibliography}{DKN2}
\bibitem {AC}
{Amrouche, C., Consiglieri, L.,}
\textit{On the stationary Oseen equations in $ \mathbb{R} ^3$,}
Comm. Math. Anal., 10 (2010), 5-29.


\bibitem{DK2}
{Deuring, P., Kra\v cmar,  S.,}
\textit{ Exterior Stationary Navier-Stokes Flows in 3D with
Non-Zero Velocity at Infinity: Approximation by Flows in Bounded
Domains,} Mathematische Nachrichten, 269-270 (2004), 86--115 .




\bibitem{DKN1}
{Deuring, P., Kra\v cmar, S., Ne\v casov\' a, \v S.,}
\textit{A representation formula for linearized stationary
incompressible viscous flows around rotating and translating
bodies},
Discrete Contin. Dyn. Syst. Ser. S 3(2)  (2010), 237--253.

\bibitem{DKN2}
{Deuring, P., Kra\v cmar, S., Ne\v casov\' a, \v S.,}
\textit{On pointwise decay of linearized stationary incompressible
viscous flow around rotating and translating bodies},
 SIAM J. Math. Anal., 43   (2011), 705--738.



\bibitem{DKN3}
{Deuring, P., Kra\v cmar, S., Ne\v casov\'a, \v S.,}
\textit{Linearized stationary incompressible flow around rotating
and translating bodies: asymptotic profile of the velocity
gradient and decay estimate of the second derivatives of the
velocity}, J. Differential Equations, 252 (2012), 459--476.


\bibitem {DKN4}
{Deuring, P., Kra\v cmar, S., Ne\v casov\'a, \v S.,}
\textit{ A linearized system describing stationary incompressible
viscous flow around rotating and translating bodies: improved
decay estimates of the velocity and its gradient}, In:
\textit{Dynamical Systems, Differential Equations and
Applications,} Vol. I, Ed. by W. Feng, Z. Feng, M. Grasselli, A.
Ibragimov, X. Lu, S. Siegmund and J. Voigt. Discrete  Contin. Dyn.
Syst. (Supplement 2011, 8th AIMS Conference, Dresden, Germany),
351-361 (2011).



\bibitem{DKN5}
{Deuring, P., Kra\v cmar, S., Ne\v casov\'a, \v S.,}
\textit{Pointwise decay of stationary rotational viscous incompressible flows with nonzero velocity at infinity.} J. Differential Equations, 255  (2013), 1576--1606.


\bibitem{DKN6}
{Deuring, P., Kra\v cmar, S., Ne\v casov\'a, \v S.,} \textit{Linearized stationary incompressible flow around rotating and translating bodies - Leray solutions}. Discrete  Contin. Dyn.
Syst., 7 (2014), 967--979.


\bibitem{Farwig2}
\newblock R. Farwig,
\newblock \emph{The stationary exterior 3D-problem of Oseen and
Navier-Stokes equations
in anisotropically weighted Sobolev spaces},
\newblock Math. Z., {211} (1992), 409--447.

\bibitem{Fa3}
{Farwig, R.,}
\textit{Estimates of lower order derivatives of viscous fluid flow past a rotating obstacle,}
Banach Center Publications, 70 (2005), 73--84.

\bibitem{FGK}
{Farwig, R., Galdi, G. P., Kyed, M.,}
\textit{Asymptotic structure of a Leray solution to the Navier-Stokes flow around a rotating body,}
Pacific J. Math., 253 (2011), 367--382.


\bibitem{FGNT}
{Farwig, R., Guenther, R. B., Ne\v casov\'a, \v S.,  Thomann,  E. A.,}
\textit{The fundamental solution of the linearized instationary Navier-Stokes
equations of motion around a rotating and translating body,} Discrete Contin. Dyn. Syst.  34  (2014), 511--529.

\bibitem{FH}
{Farwig, R., Hishida, T.,}
\textit{Stationary Navier-Stokes flow around a rotating obstacle,}
Funkcialaj Ekvacioj, 50 (2007), 371--403.

\bibitem{FH092}
{Farwig, R., Hishida, T.,}
\textit{Asymptotic profiles of steady Stokes and Navier-Stokes flows around a rotating obstacle,}
Ann. Univ. Ferrara, Sez. VII, 55 (2009), 263--277.

\bibitem{FH091}
{Farwig, R., Hishida, T.,}
\textit{Asymptotic profile of steady Stokes flow around a rotating obstacle,}
Manuscripta Math., 136 (2011), 315-338.

\bibitem{FH2}
{Farwig, R., Hishida, T.,}
\textit{Leading term at infinity of steady Navier-Stokes flow around a rotating obstacle,}
Math. Nachr., 284 (2011), 2065--2077.

\bibitem{FHM} {Farwig, R., Hishida, T., M\"uller, D.,}
\textit{ $L^q$-theory of a singular ``winding'' integral operator arising from fluid dynamics,}
Pacific J. Math., 215 (2004), 297--312.

\bibitem{FKN1}
{Farwig, R., Krbec, M., Ne\v casov\' a, \v S.,}
\textit{A weighted $L^q$ approach to Stokes flow around a rotating body,}
Ann. Univ. Ferrara, Sez. VII, 54, (2008), 61--84.

\bibitem{FKN2}
{Farwig, R., Krbec, M., Ne\v{c}asov\'a, \v{S.},}
\textit{ A weighted $L^q$-approach to Oseen flow around a rotating body,}
Math. Meth. Appl. Sci., 31 (2008), 551--574.



\bibitem{FN}
{Farwig, R., Neustupa, J.,}
\textit{On the spectrum of a Stokes-type operator arising from flow around a rotating body,}
Manuscripta Math., 122 (2007), 419--437.

\bibitem{G2}
\newblock G. P. Galdi,
\newblock ``On the motion of a rigid body in a viscous liquid: A mathematical analysis with applications,''
\newblock Handbook of Mathematical Fluid Dynamics, Volume 1, Ed. by S. Friedlander, D. Serre, Elsevier, 2002.

%



\bibitem{G3}
{Galdi, G. P.,}
\textit{Steady flow of a Navier-Stokes fluid around a rotating obstacle,}
J. Elasticity, 71 (2003), 1-31.

\bibitem{GK1}
\newblock G. P. Galdi  and M.  Kyed,
\newblock \emph{Steady-State Navier-Stokes Flows Pas a Rotating Body: Leary Solutions are Physically Reasonable},
\newblock Arch. Rat. Mech. Anal., {200} (2011), 21--58.

\bibitem{GK2}
{Galdi, G. P., Kyed, M.,}
\textit{Asymptotic behavior of a Leray solution around a rotating obstacle,}
Progress in Nonlinear Differential Equations and Their Applications, 60 (2011), 251-266.

\bibitem{Galdineu}
{Galdi, G. P.,}
\textit{An introduction to the mathematical theory of the Navier-Stokes equations. Steady-state problems (2nd ed.),}
Springer, New York e.a., 2011.

\bibitem{Ga1}
{Galdi, G.P.,}
\textit{An Introduction to the mathematical theory of
the Navier-Stokes equations. Vol. I. Linearized steady problems (rev. ed.),}
Springer, New York e.a., 1998.

\bibitem{Ga2}
{Galdi, G.P.,}
\textit{An introduction to the mathematical theory of the Navier-Stokes equations. Vol. II. Nonlinear steady problems,}
Springer, New York e.a., 1994.

\bibitem{GK3}
{Galdi, G. P., Kyed, M.,}
\textit{A simple proof of $L^q$-estimates for the steady-state Oseen and Stokes equations in a rotating frame.
Part I: strong solutions,}
 Proc. Am. Math. Soc., 141 (2013), 573-583.
\bibitem{GK3a}
{Galdi, G. P., Kyed, M.,}
\textit{A simple proof of $L^q$-estimates for the steady-state Oseen and Stokes equations in a rotating frame.
Part II: weak solutions,}
 Proc. Am. Math. Soc.,
141 (2013), 1313-1322.

\bibitem{GHH}
{Geissert, M., Heck, H., Hieber, M.,}
\textit{$L^{p}$ theory of the Navier-Stokes flow in the exterior of a moving or rotating obstacle,}
J. Reine Angew. Math., 596 (2006), 45-62.

%




\bibitem{GT}
{Guenther, R. B., Thomann, E. A.,}
\textit{The fundamental solution of the linearized Navier-Stokes equations for spinning
bodies in three spatial dimensions -- time dependent case,}
J. Math. Fluid Mech., 8 (2006), 77-98.
\bibitem{H1}
{Hishida, T.,}
\textit{An existence theorem for the Navier-Stokes flow in the exterior of a rotating obstacle,}
Arch. Rat. Mech. Anal., 150 (1999), 307-348.

\bibitem{H2}
{Hishida, T.,}
\textit{The Stokes operator with rotating effect in exterior domains,}
Analysis, 19 (1999), 51-67.

\bibitem{Hi}
{Hishida, T.,}
\textit{$L^q$ estimates of weak solutions to the stationary Stokes equations around a rotating body,}
J. Math. Soc. Japan, 58 (2006), 744-767.

\bibitem{HS}
{Hishida, T., Shibata, Y.,}
\textit{$L_p$-$L_q$ estimate of the Stokes operator and Navier-Stokes flows in the exterior of a rotating obstacle,}
RIMS K\^oky\^uroku Bessatsu, B1 (2007), 167-188.

\bibitem{KKNPS} {Kra\v cmar, S.,  Krbec, M.,  Ne\v casov\'a, \v S., Penel, P., Schumacher, K.,}
\textit{On the $L^q$-approach  with generalized anisotropic weights of the weak solution of the
Oseen flow around a rotating body,}
Nonlinear Analysis, 71 (2009), e2940-e2957.

\bibitem{KNPe1}
{Kra\v{c}mar, S., Ne\v{c}asov\'{ a}, \v{S}., Penel, P.,}
\textit{Estimates of weak solutions in anisotropically weighted Sobolev spaces to the
stationary rotating Oseen equations,}
IASME Transactions, 2 (2005), 854-861.

\bibitem{KNPe3}
{Kra\v{c}mar, S., Ne\v{c}asov\'{ a}, \v{S}., Penel, P.,}
\textit{Anisotropic $L^2$ estimates of weak solutions to the stationary Oseen type
equations in $ \mathbb{R} ^{3}$ for a rotating body,}
RIMS K\^oky\^uroku Bessatsu, B1 (2007), 219-235.

\bibitem{KNP-JAP}
{Kra\v{c}mar, S., Ne\v{c}asov\'{ a}, \v{S}., Penel, P.,}
\textit{Anisotropic $L^2$ estimates of weak solutions to the
stationary Oseen type equations in 3D -- exterior domain for a rotating body,}
J. Math. Soc. Japan, 62 (2010), 239-268.
\bibitem{KNP}
{Kra\v cmar, S., Novotn\'y, A., Pokorn\'y, M.,} \textit{Estimates of Oseen kernels in weighted $L^p $spaces}, J. Math. Soc. Japan,  53  (2001), 59-111.
\bibitem{KrPe1}
{Kra\v{c}mar, S., Penel, P.,}
\textit{Variational properties of a generic model equation in exterior 3D domains,}
Funkcialaj Ekvacioj, 47 (2004), 499-523.

\bibitem{KrPe2}
{Kra\v{c}mar, S., Penel, P.,}
\textit{New regularity results for a generic model equation in exterior 3D domains,}
Banach Center Publications Warsaw, 70 (2005), 139-155.


\bibitem{K}
{Kyed, M.,}
\textit{On the asymptotic structure of a Navier-Stokes flow past a rotating body,}
J. Math. Soc. Japan, 66 (2014), 1-16.

\bibitem{K1}
{Kyed, M.,}
\textit{Asymptotic profile of a linearized flow past a rotating body,}
 Q. Appl. Math., 71 (2013), 489-500.

\bibitem{K3}
{Kyed, M.,}
\textit{On a mapping property of the Oseen operator with rotation,}
Discrete Contin. Dynam. Syst. -- Ser. S., 6 (2013), 1315-1322.



\bibitem{Ne1}
{Ne\v{c}asov\'{ a}, \v{S}.,}
\textit{On the problem of the Stokes flow and Oseen flow in $\mathbb{R}^{3}$ with Coriolis force
arising from fluid dynamics,}
IASME Transaction, 2 (2005), 1262-1270.

\bibitem{Ne2}
{Ne\v{c}asov\'{ a}, \v{S}.,}
\textit{Asymptotic properties of the steady fall of a body in viscous fluids,}
Math. Meth. Appl. Sci., 27 (2004), 1969-1995.

\bibitem{NS}
{Ne\v{c}asov\' a, \v{S}, Schumacher, K.,}
\textit{ Strong solution to the Stokes equations of a flow around a rotating body in weighted $L^q$ spaces,}
 Math. Nachr., 13 (2011), 1701-1714.






\end{thebibliography}
\end{document}